\documentclass[12pt]{article}
\usepackage{setspace}
\usepackage{fullpage}
\usepackage{amsbsy}
\usepackage{amsmath}
\usepackage{amsfonts}
\usepackage{amsthm}
\usepackage{amssymb}
\usepackage{algorithmic}
\usepackage{algorithm}
\usepackage{enumerate}
\usepackage{graphicx}
\usepackage{multirow}
\usepackage{natbib}
\setcitestyle{round}
\usepackage{subfigure}
\usepackage{listings}
\usepackage{xcolor}
\theoremstyle{plain}
\newtheorem{theorem}{Theorem}

\newtheorem{lemma}{Lemma}
\theoremstyle{definition}
\newtheorem{definition}{Definition} 
\newtheorem{example}{Example} 
\theoremstyle{remark}
\newtheorem*{remark}{Remark}

\begin{document}

\title{\bf{Generalizing Distance Covariance to Measure and Test Multivariate Mutual Dependence}}

\author{Ze Jin,
David S. Matteson\thanks{Research support from an NSF Award (DMS-1455172), a Xerox PARC Faculty Research Award, and Cornell University Atkinson Center for a Sustainable Future (AVF-2017).}}

\date{\today}

\maketitle

\begin{abstract}

We propose three new measures of mutual dependence between multiple random vectors. Each measure is zero if and only if the random vectors are mutually independent. The first generalizes distance covariance from pairwise dependence to mutual dependence, while the other two measures are sums of squared distance covariances. The proposed measures share similar properties and asymptotic distributions with distance covariance, and capture non-linear and non-monotone mutual dependence between the random vectors. Inspired by complete and incomplete V-statistics, we define empirical and simplified empirical measures as a trade-off between the complexity and statistical power when testing mutual independence. The implementation of corresponding tests is demonstrated by both simulation results and real data examples.

\end{abstract}

{\bf Key words:}\quad {\small characteristic functions; distance covariance; multivariate analysis; mutual independence; V-statistics}

\section{Introduction}\label{intro}

Let $X = (X_1, \dots, X_d)$ be a set of variables where each component $X_j$, $j = 1, \dots, d$ is a continuous random vector, and let $\mathbf{X} = \{X^k = (X_1^k, \dots, X_d^k): k = 1, \dots, n \}$ be an i.i.d.\ sample from $F_X$, the joint distribution of $X$. We are interested in testing the hypothesis
\begin{equation*}
H_0: X_1, \dots, X_d \textrm{ are mutually independent}, \quad H_A: X_1, \dots, X_d \textrm{ are dependent},
\end{equation*}
which has many applications, including independent component analysis \citep{matteson2017independent}, graphical models \citep{fan2015conditional}, naive Bayes classifiers \citep{tibshirani2002diagnosis}, etc. This problem has been studied under different settings and assumptions, including pairwise ($d = 2$) and mutual ($d \geq 2$) independence, univariate ($X_1, \dots, X_d \in \mathbb{R}^1$) and multivariate ($X_1 \in \mathbb{R}^{p_1}, \dots, X_d \in \mathbb{R}^{p_d}$) components, and more.
Specifically, we focus on the general case that $X_1, \dots, X_d$ are not assumed jointly normal.

The most extensively studied case is pairwise independence with univariate components ($X_1, X_2 \in \mathbb{R}^1$):
Rank correlation is considered as a non-parametric counterpart to Pearson's product-moment correlation \citep{pearson1895note},
including Kendall's $\tau$ \citep{kendall1938new}, Spearman's $\rho$ \citep{spearman1904proof}, etc.
\citet{bergsma2014consistent} proposed a test based on an extension of Kendall's $\tau$, testing an equivalent condition to $H_0$.
Additionally, \citet{hoeffding1948non} proposed a non-parametric test based on marginal and joint distribution functions, testing a necessary condition to investigate $H_0$.

For pairwise independence with multivariate components ($X_1 \in \mathbb{R}^{p_1}, X_2 \in \mathbb{R}^{p_2}$):
\citet{szekely2007measuring}, \citet{szekely2009brownian} proposed a test based on distance covariance with fixed $p_1, p_2$ and $n \rightarrow \infty$, testing an equivalent condition to $H_0$. Further, \citet{szekely2013distance} proposed a $t$-test based on a modified distance covariance for the setting in which $n$ is finite and $p_1, p_2 \rightarrow \infty$, testing an equivalent condition to $H_0$ as well.

For mutual independence with univariate components ($X_1, \dots, X_d \in \mathbb{R}^1$):
One natural way to extend the pairwise rank correlation to multiple components is to collect the rank correlations between all pairs of components, and examine the norm ($\mathcal{L}_2, \mathcal{L}_\infty$) of this collection.
\citet{leung2015testing} proposed a test based on the $\mathcal{L}_2$ norm with $n, d \rightarrow \infty$, and $d/n \rightarrow \gamma \in (0, \infty)$, and \citet{han2014distribution} proposed a test based on the $\mathcal{L}_\infty$ norm with $n, d \rightarrow \infty$, and $d/n \rightarrow \gamma \in [0, \infty]$. Each are testing a necessary condition to $H_0$, in general.

For mutual independence with multivariate components ($X_1 \in \mathbb{R}^{p_1}, \dots, X_d \in \mathbb{R}^{p_d}$):
This challenging scenario has not been well studied. \citet{yao2016testing} proposed a test based on distance covariance between all pairs of components with $n, d \rightarrow \infty$, testing a necessary condition to $H_0$.
Inspired by distance covariance in \citet{szekely2007measuring}, we propose a new test based on measures of mutual dependence with fixed $d, p_1, \dots, p_d$ and $n \rightarrow \infty$ in this paper, testing an equivalent condition to $H_0$.
All computational complexities in this paper make no reference to the dimensions $d, p_1, \dots, p_d$, as they are treated as constants.

Our measures of mutual dependence involve V-statistics, and are 0 if and only if mutual independence holds.
They belong to energy statistics \citep{szekely2013energy}, and share many statistical properties with distance covariance.
Besides, \citet{pfister2016kernel} proposed $d$-variable Hilbert$-$Schmidt
independence criterion (dHSIC) under the same setting, which originates from HSIC \citep{gretton2005measuring},
and also is 0 if and only if mutual independence holds.
Although dHSIC involves V-statistics as well, they pursue kernel methods and overcome the computation bottleneck by resampling and Gamma approximation,
while we take advantage of characteristic functions and resort to incomplete V-statistics.

The weakness of testing mutual independence by a necessary condition, all pairwise independencies motivates our work on measures of mutual dependence, which is demonstrated by examples in section \ref{data}: If we directly test mutual independence based on the measures of mutual dependence proposed in this paper, we successfully detect mutual dependence. Alternatively, if we check all pairwise independencies based on distance covariance, we fail to detect any pairwise dependence, and mistakenly conclude that mutual independence holds probably because the mutual effect averages out when we narrow down to a pair.



The rest of this paper is organized as follows. In section \ref{dc}, we give a brief overview of distance covariance. In section \ref{com}, we generalize distance covariance to complete measure of mutual dependence, with its properties and asymptotic distributions derived. In section \ref{asym}, we propose asymmetric and symmetric measures of mutual dependence, defined as sums of squared distance covariances. We present synthetic and real data analysis in section \ref{data}, followed by simulation results in section \ref{sim}\footnote{An accompanying R package \texttt{EDMeasure} \citep{jin2018edm} is available on CRAN.}. Finally, section \ref{con} is the summary of our work. All proofs have been moved to appendix.

The following notations will be used throughout this paper. Let $(\cdot, \cdot, \dots, \cdot)$ denote a concatenation of (vector) components into a vector. Let $t = (t_1, \dots, t_d), t^0 = (t^0_1, \dots, t^0_d), X = (X_1, \dots, X_d) \in \mathbb{R}^{p}$ where $t_j, t^0_j, X_j \in \mathbb{R}^{p_j}$, such that $p_j$ is the marginal dimension, $j = 1, \dots, d$, and $p = \sum_{j=1}^d p_j$ is the total dimension. The assumed ``$X$'' under $H_0$ is denoted by $\widetilde{X} = (\widetilde{X}_1, \dots, \widetilde{X}_d)$, where $\widetilde{X}_j \overset{\mathcal{D}}{=} X_j$, $j = 1, \dots, d$, $\widetilde{X}_1, \dots, \widetilde{X}_d$ are mutually independent, and $X, \widetilde{X}$ are independent. Let $X', X''$ be independent copies of $X$, i.e., $X, X', X'' \overset{i.i.d.}{\sim} F_X$, and $\widetilde{X}', \widetilde{X}''$ be independent copies of $\widetilde{X}$, i.e., $\widetilde{X}, \widetilde{X}', \widetilde{X}'' \overset{i.i.d.}{\sim} F_{\widetilde{X}}$. Let the weighted $\mathcal{L}_2$ norm $\|\cdot\|_w$ of complex-valued function $\eta(t)$ be defined by $\|\eta(t)\|^2_w = \int_{\mathbb{R}^{p}} |\eta(t)|^2 w(t) \,dt$ where $|\eta(t)|^2 = \eta(t)\overline{\eta(t)}$, $\overline{\eta(t)}$ is the complex conjugate of $\eta(t)$, and $w(t)$ is any positive weight function for which the integral exists.

Given the i.i.d.\ sample $\mathbf{X}$ from $F_X$, let $\mathbf{X}_j = \{X^k_j: k = 1, \dots, n\}$ denote the corresponding i.i.d.\ sample from $F_{X_j}$, $j = 1, \dots, d$, such that $\mathbf{X} = \{\mathbf{X}_1, \dots, \mathbf{X}_d\}$. Denote the joint characteristic functions of $X$ and $\widetilde{X}$ as $\phi_{X}(t) = \textrm{E}[e^{i \langle t, X\rangle}]$ and $\phi_{\widetilde{X}}(t) = \prod_{j = 1}^d \textrm{E}[e^{i \langle t_j, X_j\rangle}]$, and denote the empirical versions of $\phi_{X}(t)$ and $\phi_{\widetilde{X}}(t)$ as $\phi^n_{X}(t) = \frac{1}{n} \sum_{k=1}^n e^{i \langle t, X^k\rangle}$ and $\phi_{\widetilde{X}}^n(t) = \prod_{j = 1}^d (\frac{1}{n} \sum_{k=1}^n e^{i \langle t_j, X^k_j\rangle})$.

\section{Distance Covariance}\label{dc}

\citet{szekely2007measuring} proposed distance covariance to capture non-linear and non-monotone pairwise dependence between two random vectors ($X_1 \in \mathbb{R}^{p_1}, X_2 \in \mathbb{R}^{p_2}$).

$X_1, X_2$ are pairwise independent if and only if $\phi_{X}(t) = \phi_{X_1}(t_1)\phi_{X_2}(t_2)$, $\forall t$, which is equivalent to $\int_{\mathbb{R}^{p}} |\phi_{X}(t) - \phi_{\widetilde{X}}(t)|^2 w(t) \,dt = 0$, $\forall w(t) > 0$ if the integral exists.
A class of the weight functions $w_0(t, m) = (K(p_1; m) K(p_2; m) |t_1|^{p_1+m} |t_2|^{p_2+m})^{-1}$ make the integral a finite and meaningful quantity composed of $m$-th moments according to Lemma 1 in \citet{szekely2005new}, where $K(q, m) = \frac{2 \pi^{q/2} \Gamma(1 - m/2) }{m 2^m \Gamma((q + m)/2)}$, and $\Gamma$ is the gamma function.

The non-negative distance covariance $\mathcal{V}(X)$ is defined by $\mathcal{V}^2(X) = \|\phi_{X}(t) - \phi_{\widetilde{X}}(t)\|^2_{w_0} = \int_{\mathbb{R}^{p}} |\phi_{X}(t) - \phi_{\widetilde{X}}(t)|^2 w_0(t) \,dt$, where
\begin{equation}
w_0(t) = (K_{p_1} K_{p_2} |t_1|^{p_1+1} |t_2|^{p_2+1})^{-1},
\end{equation}
with $m = 1$ and $K_q = K(q, 1)$, while any following result can be generalized to $0 < m < 2$.
If $\textrm{E}|X| < \infty$, then $\mathcal{V}(X) \in [0, \infty)$, and $\mathcal{V}(X) = 0$ if and only if $X_1, X_2$ are pairwise independent.

The non-negative empirical distance covariance $\mathcal{V}_n(\mathbf{X})$ is defined by $\mathcal{V}_n^2(\mathbf{X}) = \|\phi^n_{X}(t) - \phi^n_{\widetilde{X}}(t)\|^2_{w_0} = \int_{\mathbb{R}^{p}} |\phi^n_{X}(t) - \phi^n_{\widetilde{X}}(t)|^2 w_0(t) \,dt$.
Calculating $\mathcal{V}_n^2(\mathbf{X})$ via the symmetry of Euclidian distances has the time complexity $O(n^2)$.
Some asymptotic properties of $\mathcal{V}_n(\mathbf{X})$ are derived. If $\textrm{E}|X| < \infty$, then
(i) $\mathcal{V}_n(\mathbf{X}) \underset{n \rightarrow \infty}{\overset{a.s.}{\longrightarrow}} \mathcal{V}(X)$.
(ii) Under $H_0$, $n\mathcal{V}_n^2(\mathbf{X}) \underset{n \rightarrow \infty}{\overset{\mathcal{D}}{\longrightarrow}} \|\zeta(t)\|^2_{w_0}$
where $\zeta(t)$ is a complex-valued Gaussian process with mean zero and covariance function $R(t, t^0) = [\phi_{X_1}(t_1 - t_1^0) - \phi_{X_1}(t_1) \overline{\phi_{X_1}(t_1^0)}][\phi_{X_2}(t_2 - t_2^0) - \phi_{X_2}(t_2) \overline{\phi_{X_2}(t_2^0)}]$.
(iii) Under $H_A$, $n\mathcal{V}_n^2(\mathbf{X}) \underset{n \rightarrow \infty}{\overset{a.s.}{\longrightarrow}} \infty$.


\section{Complete Measure of Mutual Dependence}\label{com}

Generalizing the idea of distance covariance, we propose complete measure of mutual dependence to capture non-linear and non-monotone mutual dependence between multiple random vectors
($X_1 \in \mathbb{R}^{p_1}, \dots, X_d \in \mathbb{R}^{p_d}$).

$X_1, \dots, X_d$ are mutually independent if and only if $\phi_{X}(t) = \phi_{X_1}(t_1) \dots \phi_{X_d}(t_d) = \phi_{\widetilde{X}}(t)$, $\forall t$, which is equivalent to $\int_{\mathbb{R}^{p}} |\phi_{X}(t) - \phi_{\widetilde{X}}(t)|^2 w(t) \,dt$ $= 0$, $\forall w(t) > 0$ if the integral exists.
We put all components together instead of separating them, and choose the weight function
\begin{equation}
w_1(t) = (K_p |t|^{p+1})^{-1}.
\end{equation}

\begin{definition}\label{com-def}
The complete measure of mutual dependence $\mathcal{Q}(X)$ is defined by
\[ \mathcal{Q}(X) = \|\phi_{X}(t) - \phi_{\widetilde{X}}(t)\|^2_{w_1} = \int_{\mathbb{R}^{p}} |\phi_{X}(t) - \phi_{\widetilde{X}}(t)|^2 w_1(t) \,dt. \]
\end{definition}

We can show an equivalence to mutual independence based on $\mathcal{Q}(X)$ according to Lemma 1 in \citet{szekely2005new}.
\begin{theorem}\label{com-thm}
If $\textrm{E}|X| < \infty$, then $\mathcal{Q}(X) \in [0, \infty)$, and $\mathcal{Q}(X) = 0$ if and only if $X_1, \dots, X_d$ are mutually independent. In addition, $\mathcal{Q}(X)$ has an interpretation as expectations
\[ \mathcal{Q}(X) = \textrm{E}|X - \widetilde{X}'| + \textrm{E}|X' - \widetilde{X}| - \, \textrm{E}|X - X'| - \textrm{E}|\widetilde{X} - \widetilde{X}'|. \]
\end{theorem}

It is straightforward to estimate $\mathcal{Q}(X)$ by replacing the characteristic functions with the empirical characteristic functions from the sample.
\begin{definition}\label{emp-com-def}
The empirical complete measure of mutual dependence $\mathcal{Q}_n(\mathbf{X})$ is defined by
\[ \mathcal{Q}_n(\mathbf{X}) = \|\phi^n_{X}(t) - \phi^n_{\widetilde{X}}(t)\|^2_{w_1} = \int_{\mathbb{R}^{p}} |\phi^n_{X}(t) - \phi^n_{\widetilde{X}}(t)|^2 w_1(t) \,dt. \]
\end{definition}
\begin{lemma}\label{emp-com-thm}
$\mathcal{Q}_n(\mathbf{X})$ has an interpretation as complete V-statistics
\begin{eqnarray*}
\mathcal{Q}_n(\mathbf{X}) &=& \frac{2}{n^{d+1}}\sum_{k, \ell_1, \dots, \ell_d=1}^n |X^k - (X_1^{\ell_1}, \dots, X_d^{\ell_d})| + \, \frac{1}{n^2} \sum_{k,\ell=1}^n |X^k - X^\ell| \\
&& - \, \frac{1}{n^{2d}} \sum_{k_1, \dots, k_d, \ell_1, \dots, \ell_d = 1}^n  |(X_1^{k_1}, \dots, X_d^{k_d}) - (X_1^{\ell_1}, \dots, X_d^{\ell_d})|,
\end{eqnarray*}
whose naive implementation has the time complexity $O(n^{2d})$.
\end{lemma}

In view of the definition of distance covariance,
it may seem natural to define the measure using the weight function
\begin{equation}
w_2(t) = (K_{p_1}\dots K_{p_d} |t_1|^{p_1+1} \dots |t_d|^{p_d+1})^{-1},
\end{equation}
which equals $w_0(t)$ when $d = 2$.
Given the weight function $w_2(t)$, we can define the squared distance covariance of mutual dependence $\mathcal{U}(X) = \|\phi_{X}(t) - \phi_{\widetilde{X}}(t)\|^2_{w_2}$ and its empirical counterpart $\mathcal{U}_n(\mathbf{X}) = \|\phi^n_{X}(t) - \phi^n_{\widetilde{X}}(t)\|^2_{w_2}$,
which equal $\mathcal{V}^2(X)$ and $\mathcal{V}_n^2(\mathbf{X})$ when $d = 2$.
The naive implementation of $\mathcal{U}_n(\mathbf{X})$ has the time complexity $O(n^{d+1})$.


The reason to favor $w_1(t)$ instead of $w_2(t)$ is a trade-off between the moment condition and time complexity.
We often cannot afford the time complexity of $\mathcal{Q}_n(\mathbf{X})$ or $\mathcal{U}_n(\mathbf{X})$,
and have to simplify them through incomplete V-statistics.
An incomplete V-statistic is obtained by sampling the terms of
a complete V-statistic, where the summation extends over only a subset of the tuple of indices.
To simplify by replacing complete V-statistics with incomplete V-statistics,
$\mathcal{U}_n(\mathbf{X})$ requires the additional $d$-th moment condition $\textrm{E}|X_1  \dots  X_d| < \infty$,
while $\mathcal{Q}_n(\mathbf{X})$ does not require any other condition in addition to the first moment condition $\textrm{E}|X| < \infty$.
Thus, we can reduce the complexity of $\mathcal{Q}_n(\mathbf{X})$ to $O(n^2)$ with a weaker condition,
which makes $\mathcal{Q}(X)$ and $\mathcal{Q}_n(\mathbf{X})$ from $w_1(t)$ a more general solution.
Moreover, we define the simplified empirical version of $\phi_{\widetilde{X}}(t)$ as
\begin{equation*}
\phi_{\widetilde{X}}^{n\star}(t) = \frac{1}{n}\sum_{k=1}^n e^{i \sum_{j=1}^d \langle t_j, X_j^{k+j-1} \rangle} = \frac{1}{n}\sum_{k=1}^n e^{i \langle t, (X_1^k, \dots, X_d^{k+d-1}) \rangle},
\end{equation*}
in order to substitute $\phi_{\widetilde{X}}^{n}(t)$ for simplification, where $X_j^{n + k}$ is interpreted as $X_j^k$ for $k > 0$.
\begin{definition}\label{sim-emp-com-def}
The simplified empirical complete measure of mutual dependence $\mathcal{Q}_n^\star(\mathbf{X})$ is defined by
\[ \mathcal{Q}_n^\star(\mathbf{X}) = \|\phi_{X}^n(t) - \phi_{\widetilde{X}}^{n\star}(t) \|^2_{w_1} = \int_{\mathbb{R}^{p}} |\phi^n_{X}(t) - \phi^{n\star}_{\widetilde{X}}(t)|^2 w_1(t) \,dt. \]
\end{definition}
\begin{lemma}\label{sim-emp-com-thm}
$\mathcal{Q}^\star_n(\mathbf{X})$ has an interpretation as incomplete V-statistics
\begin{eqnarray*}
\mathcal{Q}_n^\star(\mathbf{X}) &=& \frac{2}{n^{2}}\sum_{k, \ell = 1}^n |X^k - (X_1^{\ell}, \dots, X_d^{\ell + d - 1})| + \, \frac{1}{n^2} \sum_{k,\ell=1}^n |X^k - X^\ell| \\
&& - \, \frac{1}{n^{2}} \sum_{k, \ell = 1}^n  |(X_1^{k}, \dots, X_d^{k + d - 1}) - (X_1^{\ell}, \dots, X_d^{\ell + d - 1})|,
\end{eqnarray*}
whose naive implementation has the time complexity $O(n^2)$.
\end{lemma}

Using a similar derivation to Theorem 2 and 5 of \citet{szekely2007measuring}, some asymptotic distributions of $\mathcal{Q}_n(\mathbf{X}), \mathcal{Q}^\star_n(\mathbf{X})$ are obtained as follows.
\begin{theorem}\label{emp-com-strong}
If $\textrm{E}|X| < \infty$, then
\begin{eqnarray*}
\mathcal{Q}_n(\mathbf{X}) \underset{n \rightarrow \infty}{\overset{a.s.}{\longrightarrow}} \mathcal{Q}(X) & and & \mathcal{Q}^\star_n(\mathbf{X}) \underset{n \rightarrow \infty}{\overset{a.s.}{\longrightarrow}} \mathcal{Q}(X).
\end{eqnarray*}
\end{theorem}
\begin{theorem}\label{emp-com-weak}
If $\textrm{E}|X| < \infty$, then under $H_0$, we have
\begin{eqnarray*}
n\mathcal{Q}_n(\mathbf{X}) \underset{n \rightarrow \infty}{\overset{\mathcal{D}}{\longrightarrow}} \|\zeta(t)\|^2_{w_1} & and & n\mathcal{Q}^\star_n(\mathbf{X}) \underset{n \rightarrow \infty}{\overset{\mathcal{D}}{\longrightarrow}} \|\zeta^\star(t)\|^2_{w_1},
\end{eqnarray*}
where $\zeta(t), \zeta^\star(t)$ are complex-valued Gaussian processes with mean zero and covariance functions
\begin{eqnarray*}
R(t, t^0) &=& \prod_{j = 1}^d \phi_{X_j}(t_j - t^0_j) + (d-1) \prod_{j = 1}^d \phi_{X_j}(t_j) \overline{\phi_{X_j}(t^0_j)} - \sum_{j=1}^d \phi_{X_j}(t_j - t^0_j) \prod_{\ell \neq j} \phi_{X_{\ell}}(t_{\ell}) \overline{\phi_{X_{\ell}}(t^0_{\ell})},\\
R^\star(t, t^0) &=& 2 R(t, t^0).
\end{eqnarray*}
Under $H_A$, we have
\begin{eqnarray*}
n\mathcal{Q}_n(\mathbf{X}) \underset{n \rightarrow \infty}{\overset{a.s.}{\longrightarrow}} \infty & and & n\mathcal{Q}^\star_n(\mathbf{X}) \underset{n \rightarrow \infty}{\overset{a.s.}{\longrightarrow}} \infty.
\end{eqnarray*}
\end{theorem}

Therefore, a mutual independence test can be proposed based on the weak convergence of $n\mathcal{Q}_n(\mathbf{X}), n\mathcal{Q}^\star_n(\mathbf{X})$ in Theorem \ref{emp-com-weak}. Since the asymptotic distributions of $n\mathcal{Q}_n(\mathbf{X}), n\mathcal{Q}_n^\star(\mathbf{X})$ depend on $F_X$, a permutation procedure is used to approximate them in practice.

\section{Asymmetric and Symmetric Measures of Mutual Dependence}\label{asym}

As an alternative, we now propose the asymmetric and symmetric measures of mutual dependence to capture mutual dependence via aggregating pairwise dependencies.

The subset of components on the right of $X_c$ is denoted by $X_{c^+} = (X_{c+1}, \dots, X_{d})$, with $t_{c^+} = (t_{c+1}, \dots, t_d)$, $c = 0, 1, \dots, d-1$. The subset of components except $X_c$ is denoted by $X_{-c} = ({X}_1, \dots, {X}_{c-1}, X_{c^+})$, with $t_{-c} = (t_1, \dots, t_{c-1}, t_{c^+})$, $c = 1, \dots, d-1$.

We denote pairwise independence by ${\perp\!\!\!\!\perp}$. The collection of pairwise independencies implied by mutual independence includes ``one versus others on the right''
\begin{equation}\label{onevsone}
\{X_1 {\perp\!\!\!\!\perp} X_{1^+}, X_2 {\perp\!\!\!\!\perp} X_{2^+}, \dots, X_{d-1} {\perp\!\!\!\!\perp} X_d\},
\end{equation}
``one versus all the others''
\begin{equation}\label{onevsother}
\{X_1 {\perp\!\!\!\!\perp} X_{-1}, X_2 {\perp\!\!\!\!\perp} X_{-2}, \dots, X_d {\perp\!\!\!\!\perp} X_{-d}\},
\end{equation}
and many others, e.g., $(X_1, X_2) {\perp\!\!\!\!\perp} X_{2^+}$. In fact, the number of pairwise independencies resulting from mutual independence is at least $2^{d-1} - 1$, which grows exponentially with the number of components $d$. Therefore, we cannot test mutual independence simply by checking all pairwise independencies even with moderate $d$.

Fortunately, we have two options to test only a small subset of all pairwise independencies to fulfill the task. The first one is that $H_0$ holds if and only if (\ref{onevsone}) holds, which can be verified via the sequential decomposition of distribution functions. This option is asymmetric and not unique, having $d!$ feasible subsets with respect to different orders of $X_1, \dots, X_d$. The second one is that $H_0$ holds if and only if (\ref{onevsother}) holds, which can be verified via the stepwise decomposition of distribution functions and the fact that $X_j {\perp\!\!\!\!\perp} X_{-j}$ implies $X_j {\perp\!\!\!\!\perp} X_{j^+}$. This option is symmetric and unique, having only one feasible subset.

To shed light on why these two options are necessary and sufficient conditions to mutual independence, we present the following inequality that the mutual dependence can be bounded by a sum of several pairwise dependencies as
\[ |\phi_X(t) - \prod_{j = 1}^d \phi_{X_j}(t_j)| \leq \sum_{c=1}^{d-1} |\phi_{(X_{c}, X_{c^+})}((t_{c}, t_{c^+})) - \phi_{X_{c}}(t_c) \phi_{X_{c^+}}(t_{c^+})|^2. \]

In consideration of these two options, we test a set of pairwise independencies in place of mutual independence, where we use $\mathcal{V}^2(X)$ to test pairwise independence.
\begin{definition}\label{asym-def}
The asymmetric and symmetric measures of mutual dependence $\mathcal{R}(X), \mathcal{S}(X)$ are defined by
\begin{eqnarray*}
\mathcal{R}(X) = \sum_{c=1}^{d-1} \mathcal{V}^2((X_c, X_{c^+})) & and & \mathcal{S}(X) = \sum_{c=1}^d \mathcal{V}^2((X_c, X_{-c})).
\end{eqnarray*}
\end{definition}

We can show an equivalence to mutual independence based on $\mathcal{R}(X), \mathcal{S}(X)$ according to Theorem 3 of \citet{szekely2007measuring}.
\begin{theorem}\label{asym-thm}
If $\textrm{E}|X| < \infty$, then $\mathcal{R}(X), \mathcal{S}(X) \in [0, \infty)$, and $\mathcal{R}(X), \mathcal{S}(X) = 0$ if and only if $X_1, \dots, X_d$ are mutually independent.
\end{theorem}

It is straightforward to estimate $\mathcal{R}(X), \mathcal{S}(X)$ by replacing the characteristic functions with the empirical characteristic functions from the sample.
\begin{definition}\label{emp-asym-def}
The empirical asymmetric and symmetric measures of mutual dependence $\mathcal{R}_n(\mathbf{X}), \mathcal{S}_n(\mathbf{X})$ are defined by
\begin{eqnarray*}
\mathcal{R}_n(\mathbf{X}) = \sum_{c=1}^{d-1}\mathcal{V}_n^2((\mathbf{X}_c, \mathbf{X}_{c^+})) & and & \mathcal{S}_n(\mathbf{X}) = \sum_{c=1}^d \mathcal{V}_n^2((\mathbf{X}_c, \mathbf{X}_{-c})).
\end{eqnarray*}
\end{definition}

The implementations of $\mathcal{R}_n(\mathbf{X}), \mathcal{S}_n(\mathbf{X})$ have the time complexity $O(n^2)$.
Using a similar derivation to Theorem 2 and 5 of \citet{szekely2007measuring}, some asymptotic properties of $\mathcal{R}_n(\mathbf{X}), \mathcal{S}_n(\mathbf{X})$ are obtained as follows.
\begin{theorem}\label{emp-asym-strong}
If $\textrm{E}|X| < \infty$, then
\begin{eqnarray*}
\mathcal{R}_n(\mathbf{X}) \underset{n \rightarrow \infty}{\overset{a.s.}{\longrightarrow}} \mathcal{R}(X) & and & \mathcal{S}_n(\mathbf{X}) \underset{n \rightarrow \infty}{\overset{a.s.}{\longrightarrow}} \mathcal{S}(X).
\end{eqnarray*}
\end{theorem}
\begin{theorem}\label{emp-asym-weak}
If $\textrm{E}|X| < \infty$, then
under $H_0$, we have
\begin{eqnarray*}
n\mathcal{R}_n(\mathbf{X}) \underset{n \rightarrow \infty}{\overset{\mathcal{D}}{\longrightarrow}} \sum_{j = 1}^{d-1}\|\zeta_j^R((t_j, t_{j^+}))\|^2_{w_0}& and & n\mathcal{S}_n(\mathbf{X}) \underset{n \rightarrow \infty}{\overset{\mathcal{D}}{\longrightarrow}} \sum_{j = 1}^{d}\|\zeta_j^S((t_j, t_{-j}))\|^2_{w_0},
\end{eqnarray*}
where $\zeta_j^R((t_j, t_{j^+})), \zeta_j^S((t_j, t_{-j}))$ are complex-valued Gaussian processes corresponding to the limiting distributions of $n\mathcal{V}^2_n((\mathbf{X}_j, \mathbf{X}_{j^+})), n\mathcal{V}^2_n((\mathbf{X}_j, \mathbf{X}_{-j}))$. Under $H_A$, we have
\begin{eqnarray*}
n\mathcal{R}_n(\mathbf{X}) \underset{n \rightarrow \infty}{\overset{a.s.}{\longrightarrow}} \infty & and & n\mathcal{S}_n(\mathbf{X}) \underset{n \rightarrow \infty}{\overset{a.s.}{\longrightarrow}} \infty.
\end{eqnarray*}
\end{theorem}

It is surprising to find that
$\mathcal{V}_n^2((\mathbf{X}_c, \mathbf{X}_{c^+})), c = 1, \dots, d-1$ are mutually independent asymptotically, and $\mathcal{V}_n^2((\mathbf{X}_c, \mathbf{X}_{-c})), c = 1, \dots, d$ are mutually independent asymptotically as well, which is a crucial discovery behind Theorem \ref{emp-asym-weak}.

Alternatively, we can plug in $\mathcal{Q}(X)$ instead of $\mathcal{V}^2(X)$ in Definition \ref{asym-def} and
$\mathcal{Q}_n(\mathbf{X})$ instead of $\mathcal{V}^2_n(\mathbf{X})$ in Definition \ref{emp-asym-def},
and define the asymmetric and symmetric measures $\mathcal{J}(X), \mathcal{I}(X)$ accordingly,
which equal $\mathcal{Q}(X), \mathcal{Q}_n(\mathbf{X})$ when $d = 2$.
The naive implementations of $\mathcal{J}_n(\mathbf{X}), \mathcal{I}_n(\mathbf{X})$ have the time complexity $O(n^4)$.
Similarly, we can replace $\mathcal{Q}_n(\mathbf{X})$ with $\mathcal{Q}^\star_n(\mathbf{X})$ to simplify them,
and define the simplified empirical asymmetric and symmetric measures $\mathcal{J}_n^\star(\mathbf{X}), \mathcal{I}_n^\star(\mathbf{X})$,
reducing their complexities to $O(n^2)$ without any other condition except the first moment condition $\textrm{E}|X| < \infty$.
Through the same derivations, we can show that $\mathcal{J}_n(\mathbf{X}), \mathcal{J}^\star_n(\mathbf{X})$, $\mathcal{I}_n(\mathbf{X}), \mathcal{I}^\star_n(\mathbf{X})$ have similar convergences as $\mathcal{R}_n(\mathbf{X}), \mathcal{S}_n(\mathbf{X})$ in Theorem \ref{emp-asym-strong} and \ref{emp-asym-weak}.

\section{Illustrative Examples}\label{data}

We start with two examples comparing different methods to show the value of our mutual independence tests. In practice, people usually check all pairwise dependencies to test mutual independence, due to the lack of reliable and universal mutual independence tests. It is very likely to miss the complicated mutual dependence structure, and make unsound decisions in corresponding applications assuming that mutual independence holds.

\subsection{Synthetic Data}\label{data1}

We define a triplet of random vectors $(X, Y, Z)$ on $\mathbb{R}^q \times \mathbb{R}^q \times \mathbb{R}^q$, where $X, Y \sim \mathcal{N}(0, I_q)$, $W \sim \textrm{Exp}(1/\sqrt{2})$, the first element of $Z$ is $Z_1 = \textrm{sign}(X_1Y_1)W$ and the remaining $q-1$ elements are $Z_{2:q} \sim \mathcal{N}(0, I_{q-1})$, and $X, Y, W, Z_{2:q}$ are mutually independent. Clearly, $(X, Y, Z)$ is a pairwise independent but mutually dependent triplet.

An i.i.d.\ sample of $(X, Y, Z)$ is randomly generated with sample size $n = 500$ and dimension $q = 5$. On the one hand, we test the null hypothesis $H_0: X, Y, Z$ are mutually independent using proposed measures $\mathcal{R}_n, \mathcal{S}_n, \mathcal{Q}_n^\star, \mathcal{J}_n^\star, \mathcal{I}_n^\star$. On the other hand, we test the null hypotheses $H_0^{(1)}: X {\perp\!\!\!\!\perp} Y$, $H_0^{(2)}: Y {\perp\!\!\!\!\perp} Z$, and $H_0^{(3)}: X {\perp\!\!\!\!\perp} Z$ using distance covariance $\mathcal{V}^2_n$.
An adaptive permutation size $B = 210$ is used for all tests.

As expected, mutual dependence is successfully captured, as the p-values of mutual independence tests are 0.0143 ($\mathcal{Q}_n^\star$), 0.0286 ($\mathcal{J}_n^\star$), 0 ($\mathcal{I}_n^\star$), 0.0381 ($\mathcal{R}_n$) and 0 ($\mathcal{S}_n$).
Meanwhile, the p-values of pairwise independence tests are 0.2905 ($X, Y$), 0.2619 ($Y, Z$), and 0.3048 ($X, Z$).
According to the Bonferroni correction for multiple tests among all the pairs, the significance level should be adjusted as $\alpha / 3$ for pairwise tests.
As a result, no signal of pairwise dependence is detected, and we cannot reject mutual independence.

\subsection{Financial Data}\label{data2}

\begin{figure}[!ht]
\begin{center}
\includegraphics[width=1\textwidth] {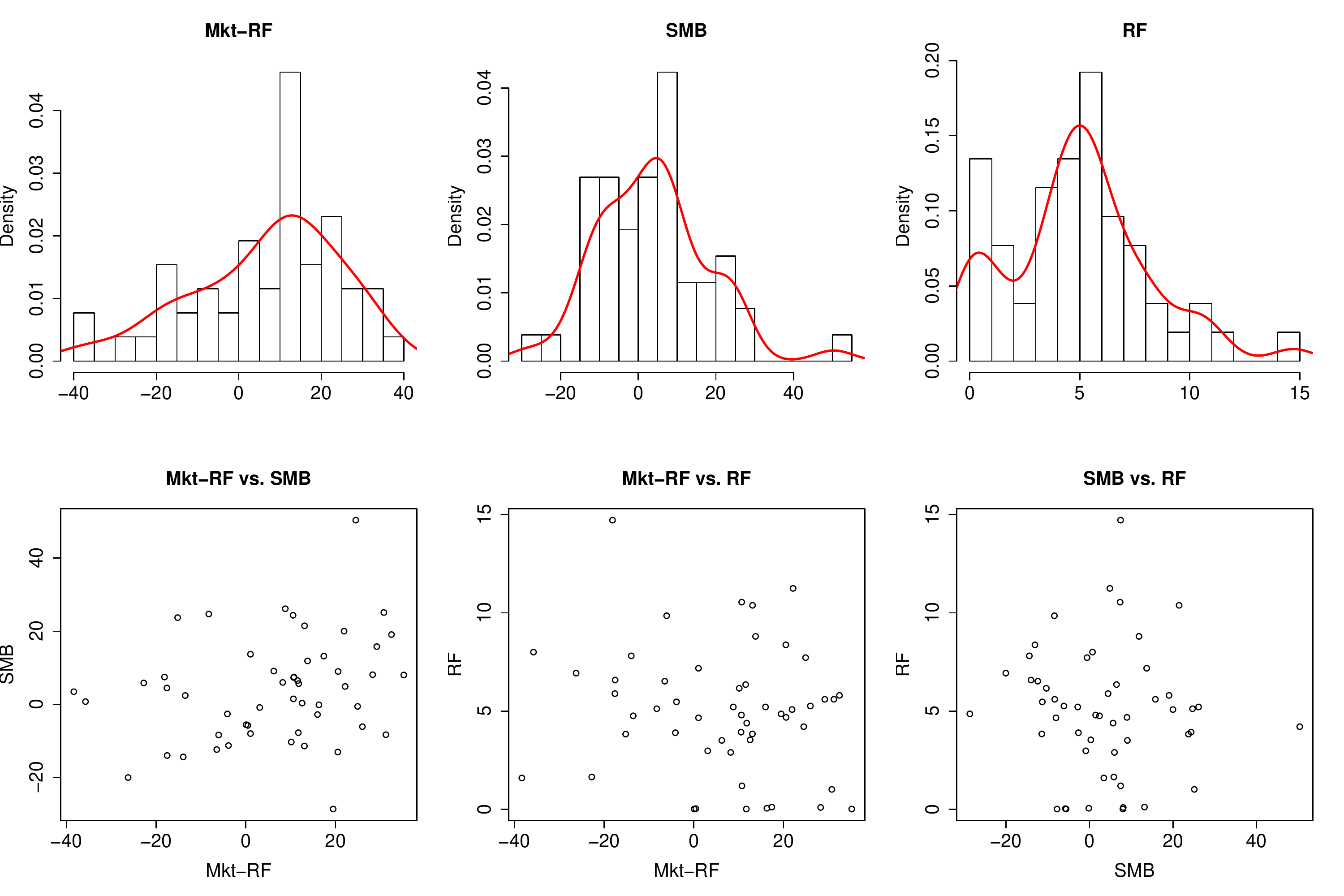}
\caption{
Three annual Fama/French factors between 1964 and 2015:
Mkt-RF (excess return on the market), SMB (small minus big) and RF (risk-free return).
The correlations are corr(Mkt-RF, SMB) = 0.238, corr(Mkt-RF, RF) = -0.161, and corr(SMB, RF) = -0.0645.
Red lines in the histograms are estimated kernel densities.
}\label{fig1}
\end{center}
\end{figure}

We collect the annual Fama/French 5 factors in the past 52 years between 1964 and 2015\footnote{Data at http://mba.tuck.dartmouth.edu/pages/faculty/ken.french/data\_library.html.}. In particular, we are interested in whether mutual dependence among three factors, $X =$ Mkt-RF (excess return on the market), $Y =$ SMB (small minus big), and $Z =$ RF (risk-free return) exists,
where annual returns are considered as nearly independent observations.
Both histograms and pair plots of $X, Y, Z$ are depicted in Figure \ref{fig1}.

For one, we apply a single mutual independence test $H_0: X, Y, Z$ are mutually independent. For another, we apply three pairwise independence tests $H_0^{(1)}: X {\perp\!\!\!\!\perp} Y$, $H_0^{(2)}: Y {\perp\!\!\!\!\perp} Z$, and $H_0^{(3)}: X {\perp\!\!\!\!\perp} Z$.
An adaptive permutation size $B = 296$ is used for all tests.

The p-values of mutual independence tests are 0.0236 ($\mathcal{Q}_n^\star$), 0.0642 ($\mathcal{J}_n^\star$), 0.0541 ($\mathcal{I}_n^\star$), 0.1588 ($\mathcal{R}_n$) and 0.1486 ($\mathcal{S}_n$), indicating that mutual dependence is successfully captured. In the meanwhile, the p-values of pairwise independence tests using distance covariance $\mathcal{V}^2_n$ are 0.1419 ($X, Y$), 0.5743 ($Y, Z$) and 0.5405 ($X, Z$).
Similarly, the significance level should be adjusted as $\alpha / 3$ according to the Bonferroni correction,
and thus we cannot reject mutual independence, since no signal of pairwise dependence is detected.

\section{Simulation Studies}\label{sim}

In this section, we evaluate the finite sample performance of proposed measures $\mathcal{Q}_n, \mathcal{R}_n, \mathcal{S}_n$, $\mathcal{J}_n, \mathcal{I}_n, \mathcal{Q}^\star_n, \mathcal{J}^\star_n, \mathcal{I}^\star_n$ by performing simulations similar to \citet{szekely2007measuring}, and compare them to benchmark measures $\mathcal{V}_n^2$ \citep{szekely2007measuring} and $\textrm{HL}^{\tau}, \textrm{HL}^{\rho}$ \citep{han2014distribution}. We also include permutation tests based on finite-sample extensions of $\textrm{HL}^{\tau}, \textrm{HL}^{\rho}$, denoted by $\textrm{HL}^{\tau}_n, \textrm{HL}^{\rho}_n$.

We test the null hypothesis $H_0$ with significance level $\alpha = 0.1$ and examine the empirical size and power of each measure. In each scenario, we run 1,000 repetitions with the adaptive permutation size $B = \lfloor 200 + 5000 / n \rfloor$ where $n$ is the sample size, for all empirical measures that require a permutation procedure to approximate their asymptotic distributions,
i.e., $\mathcal{Q}_n, \mathcal{R}_n, \mathcal{S}_n, \mathcal{J}_n, \mathcal{I}_n, \mathcal{Q}^\star_n, \mathcal{J}^\star_n, \mathcal{I}^\star_n, \mathcal{V}_n^2$.

In the following two examples, we fix $d = 2$ and change $n$ from 25 to 500,
and compare $\mathcal{Q}_n, \mathcal{R}_n, \mathcal{S}_n, \mathcal{J}_n$, $\mathcal{I}_n, \mathcal{Q}_n^\star, \mathcal{J}_n^\star, \mathcal{I}_n^\star$
to $\mathcal{V}_n^2$.

\begin{example}[pairwise multivariate normal]\label{pair-mul-nor}
$X_1, X_2 \in \mathbb{R}^5$, $(X_1, X_2)^T \sim \mathcal{N}_{10}( 0, \Sigma )$ where $\Sigma_{ii} = 1$.
Under $H_0$, $\Sigma_{ij} = 0$, $i \neq j$. Under $H_A$, $\Sigma_{ij} = 0.1$, $i \neq j$.
See results in Table \ref{t1-1} and \ref{t1-3}.
\end{example}

\begin{example}[pairwise multivariate non-normal]\label{pair-mul-non}
$X_1, X_2 \in \mathbb{R}^5$, $(Y_1, Y_2)^T \sim \mathcal{N}_{10}( 0, \Sigma )$ where $\Sigma_{ii} = 1$. $X_1 = \ln(Y_1^2), X_2 = \ln(Y_2^2)$.
Under $H_0$, $\Sigma_{ij} = 0$, $i \neq j$. Under $H_A$, $\Sigma_{ij} = 0.4$, $i \neq j$.
See results in Table \ref{t3-1} and \ref{t3-3}.
\end{example}

For both example \ref{pair-mul-nor} and \ref{pair-mul-non}, the empirical size of all measures is close to $\alpha = 0.1$.
The empirical power of $\mathcal{Q}_n, \mathcal{R}_n, \mathcal{S}_n, \mathcal{J}_n, \mathcal{I}_n$ is almost the same as that of $\mathcal{V}_n^2$, while the empirical power of $\mathcal{Q}_n^\star, \mathcal{J}_n^\star, \mathcal{I}_n^\star$ is lower than that of $\mathcal{V}_n^2$, which makes sense because we trade-off testing power and time complexity for simplified measures.

In the following two examples, we fix $d = 3$ and change $n$ from 25 to 500,
and compare $\mathcal{Q}_n, \mathcal{R}_n, \mathcal{S}_n, \mathcal{J}_n, \mathcal{I}_n$ to $\mathcal{Q}_n^\star, \mathcal{J}_n^\star, \mathcal{I}_n^\star$.

\begin{example}[mutual multivariate normal]\label{mut-mul-nor}
$X_1, X_2, X_3 \in \mathbb{R}^5$, $(X_1, X_2, X_3)^T \sim \mathcal{N}_{15}(0, \Sigma)$ where $\Sigma_{ii} = 1$.
Under $H_0$, $\Sigma_{ij} = 0$, $i \neq j$. Under $H_A$, $\Sigma_{ij} = 0.1$, $i \neq j$.
See results in Table \ref{t2-1} and \ref{t2-3}.
\end{example}

\begin{example}[mutual multivariate non-normal]\label{mut-mul-non}
$X_1, X_2, X_3 \in \mathbb{R}^5$. $(Y_1, Y_2, Y_3)^T \sim \mathcal{N}_{15}( 0, \Sigma)$ where $\Sigma_{ii} = 1$. $X_k = \ln(Y_k^2)$, $k = 1,2,3$.
Under $H_0$, $\Sigma_{ij} = 0$, $i \neq j$. Under $H_A$, $\Sigma_{ij} = 0.4$, $i \neq j$.
See results in Table \ref{t4-1} and \ref{t4-3}.
\end{example}

For both example \ref{mut-mul-nor} and \ref{mut-mul-non}, the empirical size of all measures is close to $\alpha = 0.1$.
The empirical power of $\mathcal{Q}_n, \mathcal{R}_n, \mathcal{S}_n, \mathcal{J}_n, \mathcal{I}_n$ is almost the same, the empirical power of $\mathcal{Q}_n^\star, \mathcal{J}_n^\star, \mathcal{I}_n^\star$ is almost the same, while the empirical power of $\mathcal{Q}_n^\star, \mathcal{J}_n^\star, \mathcal{I}_n^\star$ is lower than that of $\mathcal{Q}_n, \mathcal{R}_n, \mathcal{S}_n$, $\mathcal{J}_n, \mathcal{I}_n$, which makes sense since we trade-off testing power and time complexity for simplified measures.

In the last example, we change $d$ from 5 to 50 and fix $n = 100$,
and compare $\mathcal{R}_n, \mathcal{S}_n, \mathcal{Q}_n^\star$, $\mathcal{J}_n^\star, \mathcal{I}_n^\star$ to $\textrm{HL}^{\tau}, \textrm{HL}^{\rho}, \textrm{HL}^{\tau}_n, \textrm{HL}^{\rho}_n$.

\begin{example}[mutual univariate normal high-dimensional]\label{mut-uni-high}
$X_1, \dots, X_d \in \mathbb{R}^1$. $(X_1, \dots, X_d)^T \sim \mathcal{N}_{d}(0, \Sigma)$ where $\Sigma_{ii} = 1$.
Under $H_0$, $\Sigma_{ij} = 0$, $i \neq j$. Under $H_A$, $\Sigma_{ij} = 0.1$, $i \neq j$.
See results in Table \ref{t5-1} and \ref{t5-3}.
\end{example}

The empirical size of $\textrm{HL}^{\tau}, \textrm{HL}^{\rho}$ is much lower than $\alpha = 0.1$ and too conservative,
while that of other measures is fairly close to $\alpha = 0.1$.
The reason is probably that the convergence to asymptotic distributions of $\textrm{HL}^{\tau}, \textrm{HL}^{\rho}$ requires larger sample size $n$ and number of components $d$.
The measures $\mathcal{R}_n, \mathcal{S}_n$ have the highest empirical power,
and outperform the simplified measures $\mathcal{Q}_n^\star, \mathcal{J}_n^\star, \mathcal{I}_n^\star$.
The empirical power of simplified measures is similar to or even lower than that of benchmark measures when $d = 5$.
However, the empirical power of simplified measures converges much faster than that of benchmark measures as $d$ grows.

Moreover, $\mathcal{Q}_n^\star$ shows significant advantage over $\mathcal{J}_n^\star, \mathcal{I}_n^\star$. The reason is probably that $\mathcal{Q}_n^\star$ is based on truly mutual dependence while $\mathcal{J}_n^\star, \mathcal{I}_n^\star$ is based on pairwise dependencies, and large $d$ compared to $n$ introduces much more noise to $\mathcal{J}_n^\star, \mathcal{I}_n^\star$ because their summation structures, which makes them more difficult to detect mutual dependence.

The asymptotic analysis of our measures only allows small $d$ compared to $n$, while our measures work well with large $d$ compared to $n$ in example \ref{mut-uni-high}. However, this success relies on the underlying dependence structure, which is dense since each component is dependent on any other component. In contrast, if the dependence structure is sparse as each component is dependent on only a few of other components, then all measures are likely to fail.

\section{Conclusion}\label{con}

We propose three measures of mutual dependence for random vectors based on the equivalence to mutual independence through characteristic functions, following the idea of distance covariance in \citet{szekely2007measuring}.

When we select the weight function for the complete measure, we trade off between moment condition and time complexity.
Then we simplify it by replacing complete V-statistics by incomplete V-statistics, as a trade-off between testing power and time complexity.
These two trade-offs make the simplified complete measure both effective and efficient.

The asymptotic distributions of our measures depend on the underlying distribution $F_X$. Thus, the corresponding tests are not distribution-free, and we use a permutation procedure to approximate the asymptotic distributions in practice.

We illustrate the value of our measures through both synthetic and financial data examples,
where mutual independence tests based on our measures successfully capture the mutual dependence,
while the alternative checking all pairwise independencies fails and mistakenly leads to the conclusion that mutual independence holds.
Our measures achieve competitive or even better results than the benchmark measures in simulations with various examples. Although we do not allow large $d$ compared to $n$ in asymptotic analysis, our measures work well in a large $d$ example since the dependence structure is dense.

\section*{Acknowledgements}

We are grateful to Stanislav Volgushev for helpful comments on a preliminary draft of this paper.
We also thank an anonymous referee for helpful line-by-line comments.

\singlespacing


\bibliographystyle{abbrvnat}
\bibliography{Paper}

\doublespacing


\begin{table}[!ht]
\begin{center}
\caption{empirical size ($\alpha = 0.1$) in Example \ref{pair-mul-nor} with 1000 repetitions and $d = 2$.}\label{t1-1}
\begin{tabular}{|c|c|c|c|c|c|c|c|c|c|}
\hline
$n$   &$\mathcal{V}_n^2, \mathcal{R}_n, \mathcal{S}_n$ &$\mathcal{Q}_n, \mathcal{J}_n, \mathcal{I}_n$ & $\mathcal{Q}^\star_n, \mathcal{J}^\star_n$ & $\mathcal{I}^\star_n$ \\ \hline
  25 & 0.106 & 0.102 & 0.108  & 0.111\\ \hline
  30 & 0.098 & 0.115 & 0.086  & 0.114\\ \hline
  35 & 0.095 & 0.101 & 0.084  & 0.101\\ \hline
  50 & 0.101 & 0.101 & 0.111  & 0.106\\ \hline
  70 & 0.114 & 0.109 & 0.090  & 0.102\\ \hline
 100 & 0.104 & 0.105 & 0.118  & 0.117\\
\hline
\end{tabular}
\caption{empirical power ($\alpha = 0.1$) in Example \ref{pair-mul-nor} with 1000 repetitions and $d = 2$.}\label{t1-3}
\begin{tabular}{|c|c|c|c|c|c|c|c|c|c|}
\hline
$n$   &$\mathcal{V}_n^2, \mathcal{R}_n, \mathcal{S}_n$ &$\mathcal{Q}_n, \mathcal{J}_n, \mathcal{I}_n$ & $\mathcal{Q}^\star_n, \mathcal{J}^\star_n$ & $\mathcal{I}^\star_n$ \\ \hline
  25 & 0.273 & 0.246 & 0.160 &  0.182\\ \hline
  50 & 0.496 & 0.448 & 0.259 &  0.300\\ \hline
 100 & 0.807 & 0.751 & 0.442 &  0.514\\ \hline
 150 & 0.943 & 0.922 & 0.604 &  0.720\\ \hline
 200 & 0.979 &     - & 0.749 &  0.836\\ \hline
 300 & 1.000 &     - & 0.889 &  0.954\\ \hline
 500 & 1.000 &     - & 0.978 &  0.995\\
\hline
\end{tabular}
\end{center}
\end{table}

\newpage

\begin{table}[!ht]
\begin{center}
\caption{empirical size ($\alpha = 0.1$) in Example \ref{pair-mul-non} with 1000 repetitions and $d = 2$.}\label{t3-1}
\begin{tabular}{|c|c|c|c|c|c|c|c|c|c|}
\hline
$n$   &$\mathcal{V}_n^2, \mathcal{R}_n, \mathcal{S}_n$ &$\mathcal{Q}_n, \mathcal{J}_n, \mathcal{I}_n$ & $\mathcal{Q}^\star_n, \mathcal{J}^\star_n$ & $\mathcal{I}^\star_n$ \\ \hline
  25 & 0.088 & 0.093 & 0.091 & 0.092\\ \hline
  30 & 0.098 & 0.104 & 0.108 & 0.110\\ \hline
  35 & 0.104 & 0.102 & 0.104 & 0.099\\ \hline
  50 & 0.097 & 0.098 & 0.093 & 0.097\\ \hline
  70 & 0.094 & 0.097 & 0.089 & 0.097\\ \hline
 100 & 0.092 & 0.092 & 0.114 & 0.099\\
\hline
\end{tabular}
\caption{empirical power ($\alpha = 0.1$) in Example \ref{pair-mul-non} with 1000 repetitions and $d = 2$.}\label{t3-3}
\begin{tabular}{|c|c|c|c|c|c|c|c|c|c|}
\hline
$n$   &$\mathcal{V}_n^2, \mathcal{R}_n, \mathcal{S}_n$ &$\mathcal{Q}_n, \mathcal{J}_n, \mathcal{I}_n$ & $\mathcal{Q}^\star_n, \mathcal{J}^\star_n$ & $\mathcal{I}^\star_n$ \\ \hline
  25 & 0.181 & 0.185 & 0.141 & 0.152\\ \hline
  50 & 0.352 & 0.339 & 0.200 & 0.239\\ \hline
 100 & 0.610 & 0.607 & 0.372 & 0.413\\ \hline
 150 & 0.793 & 0.792 & 0.474 & 0.588\\ \hline
 200 & 0.885 &     - & 0.604 & 0.711\\ \hline
 300 & 0.989 &     - & 0.803 & 0.892\\ \hline
 500 & 0.999 &     - & 0.953 & 0.988\\
\hline
\end{tabular}
\end{center}
\end{table}


\begin{table}[!ht]
\begin{center}
\caption{empirical size ($\alpha = 0.1$) in Example \ref{mut-mul-nor} with 1000 repetitions and $d = 3$.}\label{t2-1}
\begin{tabular}{|c|c|c|c|c|c|c|c|c|}
\hline
$n$      &$\mathcal{Q}_n$ & $\mathcal{Q}_n^\star$ & $\mathcal{R}_n$ & $\mathcal{S}_n$ & $\mathcal{J}_n$ &$\mathcal{J}_n^\star$ &$\mathcal{I}_n$ &$\mathcal{I}_n^\star$ \\ \hline
  25 & 0.095 & 0.103 & 0.093 & 0.096 & 0.101 & 0.100 & 0.091 & 0.101\\ \hline
  30 &     - & 0.110 & 0.110 & 0.114 & 0.108 & 0.118 & 0.111 & 0.125\\ \hline
  35 &     - & 0.108 & 0.106 & 0.102 & 0.109 & 0.106 & 0.104 & 0.092\\ \hline
  50 &     - & 0.083 & 0.113 & 0.108 & 0.110 & 0.090 & 0.105 & 0.085\\ \hline
  70 &     - & 0.107 & 0.104 & 0.104 & 0.098 & 0.101 & 0.108 & 0.109\\ \hline
 100 &     - & 0.085 & 0.106 & 0.108 & 0.104 & 0.103 & 0.109 & 0.096\\
\hline
\end{tabular}
\caption{empirical power ($\alpha = 0.1$) in Example \ref{mut-mul-nor} with 1000 repetitions and $d = 3$.}\label{t2-3}
\begin{tabular}{|c|c|c|c|c|c|c|c|c|}
\hline
$n$      &$\mathcal{Q}_n$ & $\mathcal{Q}_n^\star$ &$\mathcal{R}_n$ & $\mathcal{S}_n$ &$\mathcal{J}_n$ &$\mathcal{J}_n^\star$ &$\mathcal{I}_n$ &$\mathcal{I}_n^\star$ \\ \hline
  25 & 0.383 & 0.220 & 0.402 & 0.418 & 0.360 & 0.199 & 0.384 & 0.228\\ \hline
  50 &     - & 0.378 & 0.707 & 0.719 & 0.651 & 0.338 & 0.671 & 0.389\\ \hline
 100 &     - & 0.707 & 0.956 & 0.961 & 0.940 & 0.643 & 0.946 & 0.767\\ \hline
 150 &     - & 0.873 & 0.996 & 0.996 & 0.993 & 0.830 & 0.994 & 0.921\\ \hline
 200 &     - & 0.946 & 1.000 & 1.000 &     - & 0.930 &     - & 0.972\\ \hline
 300 &     - & 0.997 & 1.000 & 1.000 &     - & 0.996 &     - & 0.999\\ \hline
 500 &     - & 1.000 & 1.000 & 1.000 &     - & 1.000 &     - & 1.000\\
\hline
\end{tabular}
\end{center}
\end{table}


\begin{table}[!ht]
\begin{center}
\caption{empirical size ($\alpha = 0.1$) in Example \ref{mut-mul-non} with 1000 repetitions and $d = 3$.}\label{t4-1}
\begin{tabular}{|c|c|c|c|c|c|c|c|c|}
\hline
$n$      &$\mathcal{Q}_n$ & $\mathcal{Q}_n^\star$ &$\mathcal{R}_n$ & $\mathcal{S}_n$ &$\mathcal{J}_n$ &$\mathcal{J}_n^\star$ &$\mathcal{I}_n$ &$\mathcal{I}_n^\star$ \\ \hline
  25 & 0.089 & 0.098 & 0.096 & 0.097 & 0.096 & 0.099 & 0.092 & 0.108\\ \hline
  30 &     - & 0.098 & 0.102 & 0.100 & 0.094 & 0.099 & 0.095 & 0.108\\ \hline
  35 &     - & 0.116 & 0.116 & 0.122 & 0.123 & 0.117 & 0.123 & 0.113\\ \hline
  50 &     - & 0.091 & 0.112 & 0.109 & 0.102 & 0.097 & 0.113 & 0.088\\ \hline
  70 &     - & 0.084 & 0.103 & 0.105 & 0.096 & 0.112 & 0.102 & 0.116\\ \hline
 100 &     - & 0.112 & 0.105 & 0.105 & 0.109 & 0.099 & 0.104 & 0.107\\
\hline
\end{tabular}
\caption{empirical power ($\alpha = 0.1$) in Example \ref{mut-mul-non} with 1000 repetitions and $d = 3$.}\label{t4-3}
\begin{tabular}{|c|c|c|c|c|c|c|c|c|}
\hline
$n$      &$\mathcal{Q}_n$ & $\mathcal{Q}_n^\star$ &$\mathcal{R}_n$ & $\mathcal{S}_n$ &$\mathcal{J}_n$ &$\mathcal{J}_n^\star$ &$\mathcal{I}_n$ &$\mathcal{I}_n^\star$ \\ \hline
  25 & 0.289 & 0.164 & 0.294 & 0.287 & 0.291 & 0.154 & 0.287 & 0.169\\ \hline
  50 &     - & 0.280 & 0.504 & 0.510 & 0.490 & 0.278 & 0.501 & 0.320\\ \hline
 100 &     - & 0.521 & 0.824 & 0.826 & 0.807 & 0.498 & 0.816 & 0.579\\ \hline
 150 &     - & 0.689 & 0.942 & 0.942 & 0.937 & 0.679 & 0.941 & 0.770\\ \hline
 200 &     - & 0.838 & 0.987 & 0.986 &     - & 0.826 &     - & 0.905\\ \hline
 300 &     - & 0.957 & 0.999 & 0.999 &     - & 0.956 &     - & 0.982\\ \hline
 500 &     - & 1.000 & 1.000 & 1.000 &     - & 1.000 &     - & 1.000\\
\hline
\end{tabular}
\end{center}
\end{table}

\newpage

\begin{table}[!ht]
\begin{center}
\caption{empirical size ($\alpha = 0.1$) in Example \ref{mut-uni-high} with 1000 repetitions and $n = 100$.}\label{t5-1}
\begin{tabular}{|c|c|c|c|c|c|c|c|c|c|c|c|c|}
\hline
$d$     & $\textrm{HL}^{\tau}$ & $\textrm{HL}^{\rho}$ & $\textrm{HL}^{\tau}_n$ & $\textrm{HL}^{\rho}_n$ & $\mathcal{Q}_n$ & $\mathcal{Q}_n^\star$ &$\mathcal{R}_n$ & $\mathcal{S}_n$ &$\mathcal{J}_n$ &$\mathcal{J}_n^\star$ &$\mathcal{I}_n$ & $\mathcal{I}_n^\star$ \\ \hline
  5 & 0.076 & 0.066  & 0.113  & 0.105 & - & 0.097 & 0.091 & 0.091 & - & 0.094 & - & 0.104\\ \hline
 10 & 0.077 & 0.070  & 0.104  & 0.097 & - & 0.107 & 0.092 & 0.094 & - & 0.119 & - & 0.107\\ \hline
 15 & 0.094 & 0.087  & 0.116  & 0.113 & - & 0.109 & 0.093 & 0.093 & - & 0.108 & - & 0.100\\ \hline
 20 & 0.077 & 0.066  & 0.089  & 0.089 & - & 0.096 & 0.099 & 0.118 & - & 0.115 & - & 0.101\\ \hline
 25 & 0.074 & 0.058  & 0.086  & 0.091 & - & 0.097 & 0.090 & 0.082 & - & 0.095 & - & 0.097\\ \hline
 30 & 0.091 & 0.082  & 0.110  & 0.114 & - & 0.109 & 0.092 & 0.104 & - & 0.105 & - & 0.109\\ \hline
 50 & 0.080 & 0.061  & 0.088  & 0.087 & - & 0.087 & 0.091 & 0.088 & - & 0.095 & - & 0.087\\
\hline
\end{tabular}
\caption{empirical power ($\alpha = 0.1$) in Example \ref{mut-uni-high} with 1000 repetitions and $n = 100$.}\label{t5-3}
\begin{tabular}{|c|c|c|c|c|c|c|c|c|c|c|c|c|}
\hline
$d$     & $\textrm{HL}^{\tau}$ & $\textrm{HL}^{\rho}$ & $\textrm{HL}^{\tau}_n$ & $\textrm{HL}^{\rho}_n$ & $\mathcal{Q}_n$ & $\mathcal{Q}_n^\star$ &$\mathcal{R}_n$ & $\mathcal{S}_n$ &$\mathcal{J}_n$ &$\mathcal{J}_n^\star$ &$\mathcal{I}_n$ & $\mathcal{I}_n^\star$ \\ \hline
  5 & 0.317 & 0.305  & 0.410  & 0.405 & - & 0.298 & 0.545 & 0.557 & - & 0.245 & - & 0.318\\ \hline
 10 & 0.426 & 0.416  & 0.500  & 0.510 & - & 0.557 & 0.896 & 0.915 & - & 0.409 & - & 0.497\\ \hline
 15 & 0.513 & 0.481  & 0.593  & 0.602 & - & 0.822 & 0.975 & 0.982 & - & 0.538 & - & 0.643\\ \hline
 20 & 0.558 & 0.534  & 0.625  & 0.634 & - & 0.924 & 0.996 & 0.999 & - & 0.586 & - & 0.647\\ \hline
 25 & 0.593 & 0.539  & 0.645  & 0.634 & - & 0.977 & 0.999 & 0.999 & - & 0.663 & - & 0.689\\ \hline
 30 & 0.605 & 0.556  & 0.675  & 0.664 & - & 0.980 & 1.000 & 1.000 & - & 0.711 & - & 0.700\\ \hline
 50 & 0.702 & 0.641  & 0.742  & 0.731 & - & 0.998 & 1.000 & 1.000 & - & 0.775 & - & 0.717\\
\hline
\end{tabular}
\end{center}
\end{table}

\section*{Appendix}

Proofs of Theorem \ref{com-thm}, \ref{emp-com-strong}, \ref{emp-com-weak}, \ref{asym-thm}, \ref{emp-asym-strong}, \ref{emp-asym-weak},
and Lemma \ref{emp-com-thm}, \ref{sim-emp-com-thm}.

\section*{Theorem \ref{com-thm}}

\begin{proof}

\noindent (i) $0 \leq \mathcal{Q}(X) < \infty$.

\noindent (ii) $\mathcal{Q}(X) = 0 \Longleftrightarrow X_1, \dots, X_d$ are mutually independent.

\noindent (iii) $\mathcal{Q}(X) = \textrm{E}|X - \widetilde{X}'| + \textrm{E}|X' - \widetilde{X}| - \textrm{E}|X - X'| - \textrm{E}|\widetilde{X} - \widetilde{X}'|$.

Since $w_1(t)$ is a positive weight function,
$X_1, \dots, X_d$ are mutually independent if and only if
$\mathcal{Q}(X) = \int_{\mathbb{R}^{p}} |\phi_{X}(t) - \phi_{\widetilde{X}}(t)|^2 w_1(t) \,dt$ is equal to zero.

By the boundedness property of characteristic functions and Fubini's theorem, we have

$| \phi_{X}(t) - \phi_{\widetilde{X}}(t) |^2 = \phi_{X}(t)\overline{\phi_{X}(t)} +
\phi_{\widetilde{X}}(t)\overline{\phi_{\widetilde{X}}(t)}
-\phi_{X}(t) \overline{\phi_{\widetilde{X}}(t)}
-\phi_{\widetilde{X}}(t) \overline{\phi_{X}(t)}$

$= [\textrm{E}^{i \langle t, X \rangle}] \textrm{E}[e^{-i \langle t, X \rangle}]
+ \textrm{E}[e^{i \langle t, \widetilde{X} \rangle}] \textrm{E}[e^{-i \langle t, \widetilde{X} \rangle}]
- \textrm{E}[e^{i \langle t, X \rangle}] \textrm{E}[e^{-i \langle t, \widetilde{X} \rangle}]
- \textrm{E}[e^{i \langle t, \widetilde{X} \rangle}] \textrm{E}[e^{-i \langle t, X \rangle}]$

$= \textrm{E}[e^{i\langle t, X - X' \rangle} ]
+ \textrm{E}[e^{i\langle t, \widetilde{X} - \widetilde{X}' \rangle} ]
- \textrm{E}[e^{i\langle t, X - \widetilde{X}' \rangle} ]
- \textrm{E}[e^{i\langle t, \widetilde{X} - X' \rangle} ]$

$= \textrm{E}(\cos \langle t, X - X' \rangle) + \textrm{E}(\cos \langle t, \widetilde{X} - \widetilde{X}' \rangle)
+ \textrm{E}(\cos \langle t, X - \widetilde{X}' \rangle) + \textrm{E}(\cos \langle t, \widetilde{X} - X' \rangle)$

$= \textrm{E}(1 - \cos \langle t, X - \widetilde{X}' \rangle) + \textrm{E}(1 - \cos \langle t, \widetilde{X} - X' \rangle) $

$\quad - \, \textrm{E}(1 - \cos \langle t, X - X' \rangle) - \textrm{E}(1 - \cos \langle t, \widetilde{X} - \widetilde{X}' \rangle)$.

Since $\textrm{E}(|X|) < \infty$ implies $\textrm{E}(|\widetilde{X}|) < \infty$,
we have $\textrm{E}(|X| + |\widetilde{X}|) < \infty$. Then the triangle inequality implies
$\textrm{E}|X - X'|, \textrm{E}|\widetilde{X} - \widetilde{X}'|, \textrm{E}|X - \widetilde{X}'|, \textrm{E}|\widetilde{X} - X'| < \infty$.
Therefore, by Fubini's theorem and Lemma 1, it follows that

${\cal Q}(X) = \int | \phi_{X}(t) - \phi_{\widetilde{X}}(t) |^2 \, w_1(t) \, dt$

$= \int \textrm{E}(1 - \cos \langle t, X - \widetilde{X}' \rangle) \, w_1(t) \, dt + \int \textrm{E}(1 - \cos \langle t, \widetilde{X} - X' \rangle) \, w_1(t) \, dt$

$\quad - \, \int \textrm{E}(1 - \cos \langle t, X - X' \rangle) \, w_1(t) \, dt - \int \textrm{E}(1 -  \cos \langle t, \widetilde{X} - \widetilde{X}' \rangle) \, w_1(t) \, dt$

$=\textrm{E}|X - \widetilde{X}'| + \textrm{E}|\widetilde{X} - X'| - \textrm{E}|X - X'| - \textrm{E}|\widetilde{X} - \widetilde{X}'| < \infty$.

Finally, ${\cal Q}(X) \ge 0$ since the integrand $| \phi_{X}(t) - \phi_{\widetilde{X}}(t) |^2$ is non-negative.
\end{proof}

\section*{Lemma \ref{emp-com-thm}}

\begin{proof}
After a simple calculation, we have

$|\phi_X^n(t) - \phi_{\widetilde{X}}^n(t)|^2 = \phi_X^n(t)\overline{\phi_X^n(t)} - \phi_X^n(t)\overline{\phi_{\widetilde{X}}^n(t)} - \phi_{\widetilde{X}}^n(t)\overline{\phi_X^n(t)} + \phi_{\widetilde{X}}^n(t)\overline{\phi_{\widetilde{X}}^n(t)} $

$= \frac{1}{n^2} \sum_{k, \ell = 1}^n \cos \langle t, X^k - X^\ell \rangle - \frac{2}{n^{d+1}} \sum_{k, \ell_1, \dots, \ell_d = 1}^n \cos \langle t, X^k - (X_1^{\ell_1}, \dots, X_d^{\ell_d}) \rangle$

$\quad + \, \frac{1}{n^{2d}} \sum_{k_1, \dots, k_d, \ell_1, \dots, \ell_d = 1}^n \cos \langle t, (X_1^{k_1}, \dots, X_d^{k_d}) - (X_1^{\ell_1}, \dots, X_d^{\ell_d}) \rangle + V$

$= -\frac{1}{n^2} \sum_{k, \ell = 1}^n [1 - \cos \langle t, X^k - X^\ell \rangle] + \frac{2}{n^{d+1}} \sum_{k, \ell_1, \dots, \ell_d = 1}^n [1 - \cos \langle t, X^k - (X_1^{\ell_1}, \dots, X_d^{\ell_d}) \rangle]$

$\quad - \, \frac{1}{n^{2d}} \sum_{k_1, \dots, k_d, \ell_1, \dots, \ell_d = 1}^n [1 - \cos \langle t, (X_1^{k_1}, \dots, X_d^{k_d}) - (X_1^{\ell_1}, \dots, X_d^{\ell_d}) \rangle] + V$,

\noindent where $V$ is imaginary and thus 0 as the $|\phi_X^n(t) - \phi_{\widetilde{X}}^n(t)|^2$ is real. 

By Lemma 1 in \citet{szekely2005new}

$\mathcal{Q}_n(\mathbf{X}) = \|\phi_X^n(t) - \phi_{\widetilde{X}}^n(t)\|^2_{w_1} $

$= -\frac{1}{n^2} \sum_{k, \ell = 1}^n |X^k - X^\ell| + \frac{2}{n^{d+1}} \sum_{k, \ell_1, \dots, \ell_d = 1}^n |X^k - (X_1^{\ell_1}, \dots, X_d^{\ell_d})|$

$\quad - \, \frac{1}{n^{2d}} \sum_{k_1, \dots, k_d, \ell_1, \dots, \ell_d = 1}^n |(X_1^{k_1}, \dots, X_d^{k_d}) - (X_1^{\ell_1}, \dots, X_d^{\ell_d})|$.
\end{proof}

\section*{Lemma \ref{sim-emp-com-thm}}

\begin{proof}
After a simple calculation, we have

$|\phi_X^n(t) - \phi_{\widetilde{X}}^{n\star}(t)|^2 = \phi_X^n(t)\overline{\phi_X^{n}(t)} - \phi_X^n(t)\overline{\phi_{\widetilde{X}}^{n\star}(t)} - \phi_{\widetilde{X}}^{n\star}(t)\overline{\phi_X^{n}(t)} + \phi_{\widetilde{X}}^{n\star}(t)\overline{\phi_{\widetilde{X}}^{n\star}(t)} $

$= \frac{1}{n^2} \sum_{k, \ell = 1}^n \cos \langle t, X^k - X^\ell \rangle -  \frac{2}{n^2} \sum_{k, \ell = 1}^n \cos \langle t, X^k - (X_1^{\ell}, \dots, X_d^{\ell + d - 1}) \rangle$

$\quad + \, \frac{1}{n^2} \sum_{k, \ell = 1}^n \cos \langle t, (X_1^{k}, \dots, X_d^{k + d - 1}) - (X_1^{\ell}, \dots, X_d^{\ell + d - 1}) \rangle + V^\star$

$= -\frac{1}{n^2} \sum_{k, \ell = 1}^n [1 - \cos \langle t, X^k - X^\ell \rangle] + \frac{2}{n^2} \sum_{k, \ell = 1}^n [1 - \cos \langle t, X^k - (X_1^{\ell}, \dots, X_d^{\ell + d - 1}) \rangle]$

$\quad - \, \frac{1}{n^{2}} \sum_{k, \ell = 1}^n [1 - \cos \langle t, (X_1^{k}, \dots, X_d^{k + d - 1}) - (X_1^{\ell}, \dots, X_d^{\ell + d - 1}) \rangle] + V^\star$,

\noindent where $V^\star$ is imaginary and thus 0 as the $|\phi_X^n(t) - \phi_{\widetilde{X}}^{n\star}(t)|^2$ is real. 

By Lemma Lemma 1 in \citet{szekely2005new}

$\mathcal{Q}_n^\star(\mathbf{X}) = \|\phi_X^n(t) - \phi_{\widetilde{X}}^{n\star}(t)\|^2_{w_1} $

$= -\frac{1}{n^2} \sum_{k, \ell = 1}^n |X^k - X^\ell| + \frac{2}{n^2} \sum_{k, \ell = 1}^n |X^k - (X_1^{\ell}, \dots, X_d^{\ell + d - 1})|$

$\quad - \, \frac{1}{n^{2}} \sum_{k, \ell = 1}^n |(X_1^{k}, \dots, X_d^{k + d - 1}) - (X_1^{\ell}, \dots, X_d^{\ell + d - 1})|$.
\end{proof}

\section*{Theorem \ref{emp-com-strong}}

\begin{proof}
We define
\begin{eqnarray*}
\mathcal{Q}_n = \|\phi_X^n(t) - \phi_{\widetilde{X}}^n(t)\|^2_{w_1} \triangleq \|\xi_n(t)\|^2_{w_1} & \textrm{and} & \mathcal{Q}^\star_n = \|\phi_X^n(t) - \phi_{\widetilde{X}}^{n\star}(t)\|^2_{w_1} \triangleq \|\xi_n^\star(t)\|^2_{w_1}.
\end{eqnarray*}

For $\forall 0 < \delta < 1$, define the region
\begin{eqnarray}\label{reg}
D(\delta) = \{t = (t_1, \dots, t_d): \delta \leq |t|^2 = \sum_{j=1}^d |t_j|^2 \leq 2/\delta \},
\end{eqnarray}

\noindent and random variables
\begin{eqnarray*}
\mathcal{Q}_{n, \delta} = \int_{D(\delta)} |\xi_n(t)|^2 \,dw_1 & \textrm{and} & \mathcal{Q}^\star_{n, \delta} = \int_{D(\delta)} |\xi_n^\star(t)|^2 \,dw_1.
\end{eqnarray*}

For any fixed $\delta$, the weight function $w_1(t)$ is bounded on $D(\delta)$. Hence $\mathcal{Q}_{n, \delta}$ is a combination of $V$-statistics of bounded random variables.
Similar to Theorem 2 of \citet{szekely2007measuring},
it follows by the strong law of large numbers (SLLN) for $V$-statistics \citep{mises1947asymptotic} that almost surely
\begin{eqnarray*}
\lim_{n\rightarrow \infty} \mathcal{Q}_{n, \delta} = \lim_{n\rightarrow \infty} \mathcal{Q}^\star_{n, \delta} = \mathcal{Q}_{\cdot, \delta} = \int_{D(\delta)} |\phi_X(t) - \phi_{\widetilde{X}}(t)|^2 \, dw_1.
\end{eqnarray*}

Clearly $\mathcal{Q}_{\cdot, \delta} \rightarrow \mathcal{Q}$ as $\delta \rightarrow 0$. Hence, $\mathcal{Q}_{n, \delta} \rightarrow \mathcal{Q}$ a.s.\ and $\mathcal{Q}^\star_{n, \delta} \rightarrow \mathcal{Q}$ a.s.\ as $\delta \rightarrow 0$, $n \rightarrow \infty$. In order to show $\mathcal{Q}_{n} \rightarrow \mathcal{Q}$ a.s.\ and $\mathcal{Q}^\star_{n} \rightarrow \mathcal{Q}$ a.s.\ as $n \rightarrow \infty$, it remains to prove that almost surely
\begin{eqnarray*}
\limsup_{\delta \rightarrow 0}\limsup_{n \rightarrow \infty} |\mathcal{Q}_{n, \delta} - \mathcal{Q}_{n}| = \limsup_{\delta \rightarrow 0}\limsup_{n \rightarrow \infty} |\mathcal{Q}^\star_{n, \delta} - \mathcal{Q}^\star_{n}| = 0.
\end{eqnarray*}

We define a mixture of $\widetilde{X}$ and $X$ as $Y_{-c} = (\widetilde{X}_1, \dots, \widetilde{X}_{c-1}, X_{c^+})$, $c = 1, \dots, d-1$.

By the Cauchy$-$Bunyakovsky inequality

$|\xi_n(t)|^2 = |\phi_X^n(t) - \prod_{j = 1}^d \phi_{X_j}^n(t_j)|^2$

$= |\phi_X^n(t) - \prod_{j = 1}^d \phi_{X_j}^n(t_j) - \sum_{c=1}^{d-2}  (\prod_{j=1}^c \phi_{X_j}^n(t_j) \phi_{X_{c^+}}^n(t_{c^+}) ) + \sum_{c=1}^{d-2}  (\prod_{j=1}^c \phi_{X_j}^n(t_j) \phi_{X_{c^+}}^n(t_{c^+}) ) |^2$

$\leq [|\phi_X^n(t) - \phi_{X_1}^n(t_1) \phi_{X_{1^+}}^n(t_{1^+})|$

$\quad + \, \sum_{c=1}^{d-2} |(\prod_{j=1}^c \phi_{X_j}^n(t_j)\phi_{X_{c^+}}^n(t_{c^+})) - (\prod_{j=1}^c \phi_{X_j}^n(t_j)\phi_{X_{c+1}}^n(t_{c+1}) \phi_{X_{(c+1)^+}}^n(t_{(c+1)^+}))| ]^2$

$= [\sum_{c=1}^{d-1} |\phi_{(X_{c}, Y_{-c})}^n(t_{c}, t_{-c}) - \phi_{X_{c}}^n(t_c) \phi_{Y_{-c}}^n(t_{-c})|]^2$

$\leq (d-1)\sum_{c=1}^{d-1} |\phi_{(X_{c}, Y_{-c})}^n(t_{c}, t_{-c}) - \phi_{X_{c}}^n(t_c) \phi_{Y_{-c}}^n(t_{-c})|^2$,

\noindent and

$|\xi^\star_n(t)|^2 = |\frac{1}{n} \sum_{k=1}^n e^{i \langle t, X^k \rangle} - \frac{1}{n}\sum_{k=1}^n e^{i \sum_{j=1}^d \langle t_j, X_j^{k+j-1} \rangle}|^2$

$= |\frac{1}{n} \sum_{k=1}^n ( e^{i \langle t, X^k \rangle}
- \sum_{c=2}^{d-1} e^{i \langle t, (X_1^k, \dots, X_c^{k+c-1}, X_{c^+}^k) \rangle}
+ \sum_{c=2}^{d-1} e^{i \langle t, (X_1^k, \dots, X_c^{k+c-1}, X_{c^+}^k) \rangle}
- e^{i \sum_{j=1}^d \langle t_j, X_j^{k+j-1} \rangle})|^2$

$= |\frac{1}{n} \sum_{k=1}^n
\sum_{c=1}^{d-1} (e^{i \langle t, (X_1^k, \dots, X_c^{k+c-1}, X_{c^+}^k) \rangle}
- e^{i \langle t, (X_1^k, \dots, X_{c+1}^{k+c}, X_{(c+1)^+}^k) \rangle})|^2$

$\leq (d-1) \sum_{c=1}^{d-1} |\frac{1}{n} \sum_{k=1}^n
e^{i \langle t_{-(c+1)}, (X_1^k, \dots, X_c^{k+c-1}, X_{(c+1)^+}^k) \rangle}
(e^{i \langle t_{c+1}, X_{c+1}^k \rangle} - e^{i \langle t_{c+1}, X_{c+1}^{k+c} \rangle})|^2$

$\leq (d-1) \sum_{c=1}^{d-1}
(\frac{1}{n} \sum_{k=1}^n |e^{i \langle t_{-(c+1)}, (X_1^k, \dots, X_c^{k+c-1}, X_{(c+1)^+}^k) \rangle} |^2
\frac{1}{n} \sum_{k=1}^n |e^{i \langle t_{c+1}, X_{c+1}^k \rangle} - e^{i \langle t_{c+1}, X_{c+1}^{k+c} \rangle}|^2)$

$= (d-1) \sum_{c=1}^{d-1}
(\frac{1}{n} \sum_{k=1}^n |e^{i \langle t_{c+1}, X_{c+1}^k \rangle} - e^{i \langle t_{c+1}, X_{c+1}^{k+c} \rangle}|^2)$

$\leq (d-1) \sum_{c=2}^{d} \frac{2}{n} \sum_{k=1}^n (|e^{i \langle t_{c}, X_{c}^k \rangle} - \phi_{X_{c}}(t_{c})|^2 + |\phi_{X_{c}}(t_{c}) - e^{i \langle t_{c}, X_{c}^{k+c-1} \rangle}|^2)$

$= 4(d-1)\sum_{c=2}^{d} \frac{1}{n} \sum_{k=1}^n |e^{i \langle t_{c}, X_{c}^k \rangle} - \phi_{X_{c}}(t_{c})|^2$.

By the inequality $sa + (1-s)b \geq a^s b^{1-s}$, $0 < s < 1$, $a, b > 0$, we have

$|t|^{1+p} = (|t_c|^2 + |t_{-c}|^2)^{\frac{1+p}{2}} \geq (\frac{1+p_c}{2+p}|t_c|^2 + \frac{1 + \sum_{j \neq c} p_j}{2+p}|t_{-c}|^2)^{\frac{1+p}{2}} \geq (|t_c|^{\frac{2(1+p_c)}{2+p}}  |t_c|^{\frac{2(1 + \sum_{j \neq c} p_j)}{2+p}})^{\frac{1+p}{2}}$

$= |t_c|^{\frac{1 + \sum_{j \neq c} p_j}{2+p}+p_c} |t_{-c}|^{\frac{1 + p_c}{2+p} + \sum_{j \neq c} p_j} \triangleq |t_c|^{m_c+p_c} |t_{-c}|^{m_{-c} + \sum_{j \neq c} p_j}$,

\noindent where $0 < m_c < 1$, $0 < m_{-c} < 1$ and consequently

$w_1(t) = \frac{1}{K(p,1)|t|^{1+p}} \leq \frac{K(p_c, m_c) K(\sum_{j \neq c} p_j, m_{-c})}{K(p, 1)}  \frac{1}{K(p_c, m_c)|t_c|^{m_c+p_c}} \frac{1}{K(\sum_{j \neq c} p_j, m_{-c})|t_{-c}|^{m_{-c} + \sum_{j \neq c} p_j}}$

$\triangleq C(p, p_c, \sum_{j \neq c} p_j) \frac{1}{K(p_c, m_c)|t_c|^{m_c+p_c}} \frac{1}{K(\sum_{j \neq c} p_j, m_{-c})|t_{-c}|^{m_{-c} + \sum_{j \neq c} p_j}}$,

\noindent where $C(p, p_c, \sum_{j \neq c} p_j)$ is a constant depending only on $p, p_c, \sum_{j \neq c} p_j$.

By the fact $\{ \mathbb{R}^p \backslash D(\delta) \} \subset \{|t_c|^2, |t_{-c}|^2 < \delta\} \cup \{ |t_c|^2 > 1 /\delta \} \cup \{|t_{-c}|^2 > 1 /\delta \}$ and similar steps in Theorem 2 of \citet{szekely2007measuring}, almost surely

$\limsup_{\delta \rightarrow 0}\limsup_{n \rightarrow \infty} |\mathcal{Q}_{n, \delta} - \mathcal{Q}_{n}| = \limsup_{\delta \rightarrow 0}\limsup_{n \rightarrow \infty} \int_{\mathbb{R}^p \backslash D(\delta)} |\xi_n(t)|^2 \,dw_1$

$\leq (d-1)\sum_{c=1}^{d-1} \limsup_{\delta \rightarrow 0}\limsup_{n \rightarrow \infty} \int_{\mathbb{R}^p \backslash D(\delta)}  |\phi_{(X_{c}, Y_{-c})}^n(t_{c}, t_{-c}) - \phi_{X_{c}}^n(t_c) \phi_{Y_{-c}}^n(t_{-c})|^2 \,dw_1$

$\leq C(p, p_c, \sum_{j \neq c} p_j) (d-1)\sum_{c=1}^{d-1} \limsup_{\delta \rightarrow 0}\limsup_{n \rightarrow \infty} \int_{\mathbb{R}^p \backslash D(\delta)}  |\phi_{(X_{c}, Y_{-c})}^n(t_{c}, t_{-c})$

$\quad - \, \phi_{X_{c}}^n(t_c) \phi_{Y_{-c}}^n(t_{-c})|^2 \frac{1}{K(p_c, m_c)|t_c|^{m_c+p_c}} \frac{1}{K(\sum_{j \neq c} p_j, m_{-c})|t_{-c}|^{m_{-c} + \sum_{j \neq c} p_j}} \,dt_c \,dt_{-c} $

$= 0$,

\noindent and

$\limsup_{\delta \rightarrow 0}\limsup_{n \rightarrow \infty} |\mathcal{Q}^\star_{n, \delta} - \mathcal{Q}^\star_{n}| = \limsup_{\delta \rightarrow 0}\limsup_{n \rightarrow \infty} \int_{\mathbb{R}^p \backslash D(\delta)} |\xi^\star_n(t)|^2 \,dw_1$

$\leq 4(d-1)\sum_{c=2}^{d} \limsup_{\delta \rightarrow 0} \limsup_{n \rightarrow \infty} \frac{1}{n} \sum_{k=1}^n  \int_{\mathbb{R}^p \backslash D(\delta)}   |e^{i \langle t_{c}, X_{c}^k \rangle} - \phi_{X_{c}}(t_{c})|^2 \,dw_1$

$\leq C(p, p_c, \sum_{j \neq c} p_j) 4(d-1)\sum_{c=2}^{d} \limsup_{\delta \rightarrow 0}\limsup_{n \rightarrow \infty} \frac{1}{n} \sum_{k=1}^n \int_{\mathbb{R}^p \backslash D(\delta)}   |e^{i \langle t_{c}, X_{c}^k \rangle} - \phi_{X_{c}}(t_{c})|^2 $

$\, \frac{1}{K(p_c, m_c)|t_c|^{m_c+p_c}} \frac{1}{K(\sum_{j \neq c} p_j, m_{-c})|t_{-c}|^{m_{-c} + \sum_{j \neq c} p_j}} \,dt_c \,dt_{-c} $

$= 0$.

Therefore, almost surely

$\limsup_{\delta \rightarrow 0}\limsup_{n \rightarrow \infty} |\mathcal{Q}_{n, \delta} - \mathcal{Q}_{n}| = \limsup_{\delta \rightarrow 0}\limsup_{n \rightarrow \infty} |\mathcal{Q}^\star_{n, \delta} - \mathcal{Q}^\star_{n}| = 0$.
\end{proof}

\section*{Theorem \ref{emp-com-weak}}

\begin{proof}
(i) Under $H_0$:

Let $\zeta(t)$ denote a complex-valued Gaussian processe with mean zero and covariance functions
\begin{eqnarray*}
R(t, t^0) &=& \prod_{j = 1}^d \phi_{X_j}(t_j - t^0_j) + (d-1) \prod_{j = 1}^d \phi_{X_j}(t_j) \overline{\phi_{X_j}(t^0_j)}\\
&& - \, \sum_{j=1}^d \phi_{X_j}(t_j - t^0_j) \prod_{j' \neq j} \phi_{X_{j'}}(t_{j'}) \overline{\phi_{X_{j'}}(t^0_{j'})}.
\end{eqnarray*}

We define
\begin{equation*}
n \mathcal{Q}_n = n\|\phi_X^n(t) - \phi_{\widetilde{X}}^{n}(t)\|^2_{w_1} \triangleq \|\zeta_n(t) \|^2_{w_1}.
\end{equation*}

After a simple calculation, we have

$\textrm{E}[\zeta_n(t)] = \textrm{E}[\zeta_n^\star(t)] = 0$,

$\textrm{E}[\zeta_n(t)\overline{\zeta_n(t^0)}]$

$= (1- \frac{1}{n^{d-1}})\prod_{j=1}^d \phi_{X_j}(t_j-t_j^0) + (n-1- \frac{(n-1)^d}{n^{d-1}}) \prod_{j=1}^d \phi_{X_j}(t_j) \overline{\phi_{X_j}(t^0_j)}$

$\quad - \, \frac{(n-1)^{d-1}}{n^{d-1}} [\sum_{j=1}^d  \phi_{X_j}(t_j -t_j^0)  \prod_{j' \neq j}  \phi_{X_j}(t_j) \overline{\phi_{X_j}(t^0_j)}] + o_n(1)  $

$\rightarrow R(t, t^0)$ as $n \rightarrow \infty$.

In particular, $\textrm{E}|\zeta_n(t)|^2 \rightarrow R(t, t) \leq d$ as $n \rightarrow \infty$. Thus, $\textrm{E}|\zeta_n(t)|^2 \leq d + 1$ for enough large $n$.

For $\forall 0 < \delta < 1$, define the region $D(\delta)$ as (\ref{reg}). Given $\forall \epsilon > 0$, we choose a partition $\{D^\ell(\delta)\}_{\ell=1}^N$ of $D(\delta)$ into $N(\epsilon)$ measurable sets with diameter at most $\epsilon$, and suppress the notation of $D(\delta), D^\ell(\delta)$ as $D, D^\ell$.
Then we define two sequences of random variables for any fixed $t^\ell \in D^\ell, \ell = 1, \dots, N$
\begin{equation*}
Q_n(\delta) = \sum_{\ell=1}^N \int_{D^\ell} |\zeta_n(t^\ell)|^2 \,dw_1.
\end{equation*}

For any fixed $M > 0$, let $\beta(\epsilon) = \sup_{t, t^0} E||\zeta_n(t) |^2 - |\zeta_n(t^0) |^2|$ where the supremum is taken over all $t = (t_1, \dots, t_d)$ and $t^0 = (t^0_1, \dots, t^0_d)$ s.t. $\max \{|t|^2, |t^0|^2\} \leq M$ and $|t-t^0|^2 = \sum_{j=1}^d |t_j - t^0_j|^2 \leq \epsilon^2$.
By the continuous mapping theorem and $\zeta_n(t) \rightarrow \zeta_n(t^0)$ as $\epsilon \rightarrow 0$, we have
$|\zeta_n(t)|^2 \rightarrow |\zeta_n(t^0)|^2$ as $\epsilon \rightarrow 0$.
By the dominated convergence theorem and $\textrm{E}|\zeta_n(t)|^2 \leq d+1$ for enough large $n$, we have
$E||\zeta_n(t)|^2 - |\zeta_n(t^0)|^2| \rightarrow 0$ as $\epsilon \rightarrow 0$,
which leads to $\beta(\epsilon) \rightarrow 0$ as $\epsilon \rightarrow 0$.
%
%
%
%
%
%

As a result

$\textrm{E}|\int_D |\zeta_n(t)|^2 \,dw_1 - Q_n(\delta)| = \textrm{E}|\sum_{\ell=1}^N \int_{D^\ell} (|\zeta_n(t)|^2 - |\zeta_n(t^\ell)|^2) \,dw_1|$

$\leq \sum_{\ell=1}^N \int_{D^\ell} \textrm{E}||\zeta_n(t)|^2 - |\zeta_n(t^\ell)|^2| \,dw_1 \leq \beta(\epsilon) \int_D 1 \,dw_1$

$\rightarrow 0$ as $\epsilon \rightarrow 0$.

By similar steps in Theorem \ref{emp-com-strong}, we have

$\textrm{E}|\int_D |\zeta_n(t)|^2 \,dw_1 - \|\zeta_n\|^2_{w_1}| \rightarrow 0$ as $\delta \rightarrow 0$ and $\textrm{E}|\int_D |\zeta_n^\star(t)|^2 \,dw_1 - \|\zeta^\star_n\|^2_{w_1}| \rightarrow 0$ as $\delta \rightarrow 0$.

Therefore

$\textrm{E}|Q_n(\delta) - \|\zeta_n\|^2_{w_1}| \rightarrow 0$ as $\epsilon, \delta \rightarrow 0$ and $\textrm{E}|Q_n^\star(\delta) - \|\zeta^\star_n\|^2_{w_1}| \rightarrow 0$ as $\epsilon, \delta \rightarrow 0$.

On the other hand, define two random variables for any fixed $t^\ell \in D^\ell$, $\ell = 1, \dots, N$
\begin{equation*}
Q(\delta) = \sum_{\ell = 1}^N \int_{D^\ell} |\zeta(t^\ell)|^2 \,dw_1.
\end{equation*}

Similarly, we have

$\textrm{E}|Q(\delta) - \|\zeta\|^2_{w_1}| \rightarrow 0$ as $\epsilon, \delta \rightarrow 0$.

By the multivariate central limit theorem, delta method and continuous mapping theorem, we have

$Q_n(\delta) = \sum_{\ell=1}^N \int_{D^\ell} |\zeta_n(t^\ell)|^2 \,dw_1 \rightarrow_\mathcal{D} \sum_{\ell=1}^N \int_{D^\ell} |\zeta(t^\ell)|^2 \,dw_1 = Q(\delta)$ as $n \rightarrow \infty$.

Therefore

$\|\zeta_n\|^2_{w_1} \rightarrow_\mathcal{D} \|\zeta\|^2_{w_1}$ as $\epsilon, \delta \rightarrow 0$, $n \rightarrow \infty$,

\noindent since $\{Q_n(\delta)\}$ have the following properties

(a) $Q_n(\delta)$ converges in distribution to $Q(\delta)$ as $n \rightarrow \infty$.

(b) $\textrm{E}|Q_n(\delta) - \|\zeta_n\|^2_{w_1}| \rightarrow 0$ as $\epsilon, \delta \rightarrow 0$.

(c) $\textrm{E}|Q(\delta) - \|\zeta\|^2_{w_1}| \rightarrow 0$ as $\epsilon, \delta \rightarrow 0$.

Analogous to $\zeta(t), \zeta_n(t), \beta(\epsilon), Q(\delta)$, $Q_n(\delta)$ for $\mathcal{Q}_n$, we can define $\zeta^\star(t), \zeta^\star_n(t), \beta^\star(\epsilon), Q^\star(\delta)$, $Q_n^\star(\delta)$ for $\mathcal{Q}^\star_n$,
and prove that $\|\zeta^\star_n\|^2_{w_1} \rightarrow_\mathcal{D} \|\zeta^\star\|^2_{w_1}$ as $\epsilon, \delta \rightarrow 0$, $n \rightarrow \infty$
through the same derivations.
The only differences are $\textrm{E}[\zeta^\star_n(t)\overline{\zeta^\star_n(t^0)}] = 2R(t, t^0)$ and $\textrm{E}|\zeta^\star_n(t)|^2 = 2R(t, t) \leq 2d + 1$ for enough large $n$.

\noindent (ii) Under $H_A$:

By Theorem \ref{com-thm} and \ref{emp-com-strong}, we have

$\mathcal{Q}_{n} \rightarrow \mathcal{Q} > 0$ a.s.\ as $n \rightarrow \infty$.

Therefore

$n\mathcal{Q}_{n} \rightarrow \infty$ a.s.\ as $n \rightarrow \infty$.

Similarly, we can prove that $n\mathcal{Q}^\star_{n} \rightarrow \infty$ a.s.\ as $n \rightarrow \infty$
through the same derivations.
\end{proof}

\section*{Theorem \ref{asym-thm}}

\begin{proof}
\noindent (i) $0 \leq \mathcal{R}(X) < \infty$.

\noindent (ii) $0 \leq \mathcal{S}(X) < \infty$.

\noindent (iii) $\mathcal{R}(X) = \sum_{c=1}^{d-1} \mathcal{V}^2(X_c, X_{c^+}) = 0 \Longleftrightarrow X_1, \dots, X_d$ are mutually independent.

\noindent (iv) $\mathcal{S}(X) = \sum_{c=1}^{d} \mathcal{V}^2(X_c, X_{-c}) = 0 \Longleftrightarrow X_1, \dots, X_d$ are mutually independent.

Since $\textrm{E}|X| < \infty$, we have $0 \leq \mathcal{V}^2(X_c, X_{c^+}) < \infty$, $c = 1, \dots, d-1$.
Thus, $0 \leq \mathcal{R}(X) = \sum_{c=1}^{d-1} \mathcal{V}^2(X_c, X_{c^+}) < \infty$.

Similarly, we have $0 \leq \mathcal{S}(X) = \sum_{c=1}^{d} \mathcal{V}^2(X_c, X_{-c}) < \infty$.

``$\Longleftarrow$''

If $X_1, \dots, X_d$ are mutually independent, then $X_c$ and $X_{c^+}$ are independent, $\forall c = 1, \dots, d-1$.

By Theorem 3 of \citet{szekely2007measuring}, $\mathcal{V}^2(X_c, X_{c^+}) = 0$, $\forall c = 1, \dots, d-1$.

As a result, $\mathcal{R}(X) = 0$.

Similarly, we can prove that $\mathcal{S}(X) = 0$, since $X_c$ and $X_{-c}$ are independent, $\forall c = 1, \dots, d$.

``$\Longrightarrow$''

If $\mathcal{R}(X) = 0$, then $\mathcal{V}^2(X_c, X_{c^+}) = 0$, $\forall c = 1, \dots, d-1$.

By Theorem 3 of \citet{szekely2007measuring}, $X_c$ and $X_{c^+}$ are independent,
$\forall c = 1, \dots, d-1$.
Thus, For all $t \in \mathbb{R}^p$, we have
\begin{equation*}
\phi_{(X_j,X_{j^+})}(t_j,t_{j^+})  - \phi_{X_j}(X_j) \phi_{X_{j^+}}(t_{j^+}) | = 0,
\end{equation*}
where $\phi_{X_j}$ and $\phi_{X_{j+}}$ denote the marginal and
$\phi_{(X_j, X_{j^+})}$ denotes the joint characteristic function of $X_j$ and $X_{j^+}$ respectively,
$j = 1, \dots, d$.

For all $t \in \mathbb{R}^p$, we have

$| \phi_{X}(t) - \prod_{j=1}^d \phi_{X_j}(t_j)|$

$= |\phi_X(t) - \prod_{j = 1}^d \phi_{X_j}(t_j) - \sum_{c=1}^{d-2}  (\prod_{j=1}^c \phi_{X_j}(t_j) \phi_{X_{c^+}}(t_{c^+}) ) + \sum_{c=1}^{d-2}  (\prod_{j=1}^c \phi_{X_j}(t_j) \phi_{X_{c^+}}(t_{c^+}) ) |$

$\leq |\phi_X(t) - \phi_{X_1}(t_1) \phi_{X_{1^+}}(t_{1^+})|$

$\quad + \, \sum_{c=1}^{d-2} |\prod_{j=1}^c \phi_{X_j}(t_j)\phi_{X_{c^+}}(t_{c^+}) - \prod_{j=1}^c \phi_{X_j}(t_j)\phi_{X_{c+1}}(t_{c+1}) \phi_{X_{(c+1)^+}}(t_{(c+1)^+})| $

$\leq |\phi_X(t) - \phi_{X_1}(t_1) \phi_{X_{1^+}}(t_{1^+})|$

$\quad + \, \sum_{c=1}^{d-2} |\prod_{j=1}^c \phi_{X_j}(t_j)|  |\phi_{X_{c^+}}(t_{c^+}) - \phi_{X_{c+1}}(t_{c+1}) \phi_{X_{(c+1)^+}}(t_{(c+1)^+})| $

$\leq |\phi_X(t) - \phi_{X_1}(t_1) \phi_{X_{1^+}}(t_{1^+})| + \sum_{c=1}^{d-2} |\phi_{X_{c^+}}(t_{c^+}) - \phi_{X_{c+1}}(t_{c+1}) \phi_{X_{(c+1)^+}}(t_{(c+1)^+})| $

$= \sum_{c = 1}^{d-1} | \phi_{(X_c,X_{c^+})}(t_c,t_{c^+})  - \phi_{X_c}(t_c) \phi_{X_{c^+}}(t_{c^+}) |$

$= 0$.

Therefore, for all $t \in \mathbb{R}^p$, we have $|\phi_{X}(t) - \prod_{j=1}^d \phi_{X_j}(t_j)| = 0$,
which implies that $X_1, \dots, X_d$ are mutually independent.

Similarly, we can prove that $\mathcal{S}(X) = 0$ implies that $X_1, \dots, X_d$ are mutually independent,
since $X_c$ and $X_{-c}$ are independent implies that $X_c$ and $X_{c^+}$ are independent.
\end{proof}

\section*{Theorem \ref{emp-asym-strong}}

\begin{proof} By Theorem 2 of \citet{szekely2007measuring}

$\lim_{n \rightarrow \infty} \mathcal{V}_n^2(\mathbf{X_c}, \mathbf{X_{c^+}}) = \mathcal{V}^2({X_c}, {X_{c^+}})$, $c = 1, \dots, d - 1$,

$\lim_{n \rightarrow \infty} \mathcal{V}_n^2(\mathbf{X_c}, \mathbf{X_{-c}}) = \mathcal{V}^2({X_c}, {X_{-c}})$, $c = 1, \dots, d$.

Therefore, the limit of sum converges to the sum of limit as

$\mathcal{R}_n(\mathbf{X}) \underset{n \rightarrow \infty}{\overset{a.s.}{\longrightarrow}} \mathcal{R}(X)$  and  $\mathcal{S}_n(\mathbf{X}) \underset{n \rightarrow \infty}{\overset{a.s.}{\longrightarrow}} \mathcal{S}(X)$.
\end{proof}

\section*{Theorem \ref{emp-asym-weak}}

\begin{proof}
\noindent (i) Under $H_0$:

We define
\begin{eqnarray*}
n \mathcal{R}_n(\mathbf{X}) = n\sum_{c=1}^{d-1} \mathcal{V}^2_n(\mathbf{X}_c, \mathbf{X}_{c^+}) \triangleq \sum_{c=1}^{d-1}\|\zeta_n^c(t_{(c-1)^+}) \|^2_{w_0},
\end{eqnarray*}

\noindent which is the sum corresponding to the pairs $\{ X_{d-1}, X_d \}$, $\{ X_{d-2}, (X_{d-1}, X_d) \}$,
$\{ X_{d-3}, (X_{d-2},$ $X_{d-1}, X_d) \}$, $\dots$, $\{ X_1, (X_2, \dots, X_d) \}$. Any two of them can be reorganized as $\{ X_1, X_2 \}$ and $\{ X_4, (X_1, X_2, X_3) \}$ where $X_3$ could be empty. Without loss of generality, next we will show $\phi^n_{(X_1, X_2)}(t_1, t_2) - \phi_{X_1}^n(t_1) \phi_{X_2}^n(t_2)$ and $\phi^n_{(X_1, X_2, X_3, X_4)}(s_1, s_2) - \phi_{(X_1, X_2, X_3)}^n(s_1) \phi_{X_4}^n(s_2)$ are uncorrelated. Then it follows that $\zeta_n^c(t_{(c-1)^+})$, $c = 1, \dots, d-1$ are uncorrelated.

After a simple calculation, we have

$\textrm{E} [\phi^n_{(X_1, X_2)}(t_1, t_2) - \phi_{X_1}^n(t_1) \phi_{X_2}^n(t_2)] = \textrm{E} [\phi^n_{(X_1, X_2, X_3, X_4)}(s_1, s_2) - \phi_{(X_1, X_2, X_3)}^n(s_1) \phi_{X_4}^n(s_2)] = 0$,

$\textrm{E} [\phi^n_{(X_1, X_2)}(t_1, t_2) - \phi_{X_1}^n(t_1) \phi_{X_2}^n(t_2)]
[\overline{\phi^n_{(X_1, X_2, X_3, X_4)}(s_1, s_2) - \phi_{(X_1, X_2, X_3)}^n(s_1) \phi_{X_4}^n(s_2)}] = 0$.

As a result

$\textrm{Cov}(\phi^n_{(X_1, X_2)}(t_1, t_2) - \phi_{X_1}^n(t_1) \phi_{X_2}^n(t_2),
\overline{\phi^n_{(X_1, X_2, X_3, X_4)}(s_1, s_2) - \phi_{(X_1, X_2, X_3)}^n(s_1) \phi_{X_4}^n(s_2) }) = 0$.

For $\forall \delta > 0$, define the region $D_c(\delta) = \{t_{(c-1)^+} = (t_c, t_{c^+}) = (t_c, \dots, t_d): \delta \leq |t_{(c-1)^+}|^2 = \sum_{j=c}^d |t_j|^2 \leq 2/\delta \}$. Given $\forall \epsilon > 0$, we choose a partition $\{D^\ell_c \}_{\ell=1}^{N_c}$ of $D_c(\delta)$ into $N_c(\epsilon)$ measurable sets with diameter at most $\epsilon$, and define a sequence of random variables for any fixed $t^\ell_{(c-1)^+} \in D^\ell_c$, $\ell = 1, \dots, N_c$ as
\begin{eqnarray*}
Q_n^c(\delta) = \sum_{\ell=1}^{N_c} \int_{D^\ell_c} |\zeta_n^c(t^\ell_{(c-1)^+} )|^2 \,dw_0.
\end{eqnarray*}

Let $\zeta^c(t_{(c-1)^+}) = \zeta^c(t_c, t_{c^+})$ denote a complex-valued Gaussian process with mean zero and covariance function $R_c^\zeta(t_{(c-1)^+}, t^0_{(c-1)^+}) = [\phi_{X_c}(t_{c} - t^0_{c}) - \phi_{X_c}(t_c) \overline{\phi_{X_{c}}(t^0_c)}][\phi_{X_{c^+}}(t_{c^+} - t^0_{c^+}) - \phi_{X_{c^+}}(t_{c^+}) \overline{\phi_{X_{c^+}}(t^0_{c^+})}]$.

By the multivariate central limit theorem, delta method and continuous mapping theorem, we have

$\begin{pmatrix}
  Q_n^1(\delta) - \sum_{\ell=1}^{N_1} \int_{D^\ell_1} |\zeta^1(t^\ell)|^2 \,dw_0\\
  \vdots \\
  Q_n^{d-1}(\delta) - \sum_{\ell=1}^{N_{d-1}} \int_{D^\ell_{d-1}} |\zeta^{d-1}(t^\ell_{(d-2)^+})|^2 \,dw_0 \\
\end{pmatrix} \rightarrow_\mathcal{D} \begin{pmatrix}
  \sum_{\ell=1}^{N_1} \int_{D^\ell_1} |\zeta^1(t^\ell)|^2  \,dw_0\\
  \vdots \\
\sum_{\ell=1}^{N_{d-1}} \int_{D^\ell_{d-1}} |\zeta^{d-1}(t^\ell_{(d-2)^+})|^2  \,dw_0\\
\end{pmatrix}$

\noindent as $n \rightarrow \infty$ with asymptotic mutual independence.

Thus, $Q_n^c(\delta)$, $c = 1, \dots, d-1$ are asymptotically mutually independent.

By similar steps in Theorem 5 of \citet{szekely2007measuring}, we have

$\textrm{E} | Q_n^c(\delta) - \|\zeta_n^c(t_{(c-1)^+}) \|^2_{w_0}| \rightarrow 0$, $c = 1, \dots, d-1$ as $\epsilon, \delta \rightarrow 0$.

Hence

$\begin{pmatrix}
  \|\zeta_n^1(t) \|^2_{w_0} - Q_n^1(\delta) \\
  \vdots \\
\|\zeta_n^{d-1}(t_{(d-2)^+}) \|^2_{w_0} - Q_n^{d-1}(\delta) \\
\end{pmatrix}
\rightarrow_{\mathcal{P}} \begin{pmatrix}
  0\\
  \vdots \\
0\\
\end{pmatrix}$ as $\epsilon, \delta \rightarrow 0$.

By the multivariate Slutsky's theorem, we have

$\begin{pmatrix}
  \|\zeta_n^1(t) \|^2_{w_0} \\
  \vdots \\
\|\zeta_n^{d-1}(t_{(d-2)^+}) \|^2_{w_0}\\
\end{pmatrix} \rightarrow_{\mathcal{D}} \begin{pmatrix}
  \|\zeta^1(t)\|^2_{w_0} \\
  \vdots \\
\|\zeta^{d-1}(t_{(d-2)^+})\|^2_{w_0}\\
\end{pmatrix}$

\noindent as $\epsilon, \delta \rightarrow 0$, $n \rightarrow \infty$ with asymptotic mutual independence.

Therefore

$\|\zeta_n^c(t_{(c-1)^+}) \|^2_{w_0}$, $c = 1, \dots, d-1$ are asymptotically mutually independent.

$\sum_{c=1}^{d-1}\|\zeta_n^c(t_{(c-1)^+}) \|^2_{w_0} \rightarrow_{\mathcal{D}} \sum_{c=1}^{d-1}\|\zeta^c(t_{(c-1)^+})\|^2_{w_0}$ as $n \rightarrow \infty$.

Analogous to $\zeta_n^c(t_{(c-1)^+}), \zeta^c(t_{(c-1)^+}), R_c^\zeta(t_{(c-1)^+}, t^0_{(c-1)^+})$ for $\mathcal{R}_n(\mathbf{X})$, we can define $\eta_n^c(t)$, $\eta^c(t), R_c^\eta(t, t^0)$ for $\mathcal{S}_n(\mathbf{X})$, and
prove that $\|\eta_n^c(t)\|^2_{w_0}$, $c = 1, \dots, d$ are asymptotically mutually independent,
and $\sum_{c=1}^{d}\|\eta_n^c(t) \|^2_{w_0} \rightarrow_{\mathcal{D}} \sum_{c=1}^{d}\|\eta_c(t)\|^2_{w_0}$ as $n \rightarrow \infty$
through the same derivations.

The only differences are that
we will show $\phi^n_{(X_1, X_2, X_3)}(t_1, t_2, t_3) - \phi_{X_1}^n(t_1) \phi_{(X_2, X_3)}^n(t_2, t_3)$ and $\phi^n_{(X_1, X_2, X_3)}(s_1, s_2, s_3) - \phi_{X_2}^n(s_2) \phi_{(X_1, X_3)}^n(s_1, s_3)$ are asymptotically uncorrelated.

\noindent (ii) Under $H_A$:

By Theorem \ref{asym-thm}, we have

$\mathcal{R}_{n} \rightarrow \mathcal{R} > 0$ a.s.\ as $n \rightarrow \infty$.

Therefore

$n\mathcal{R}_{n} \rightarrow \infty$ a.s.\ as $n \rightarrow \infty$.

Similarly, we can prove that $n\mathcal{S}_{n} \rightarrow \infty$ a.s.\ as $n \rightarrow \infty$ through the same derivations.
\end{proof}

\begin{remark}
Under $H_A$, $\zeta_n^c(t_{(c-1)^+})$, $c = 1, \dots, d-1$ are not asymptotically uncorrelated, and $\eta_n^c(t)$, $c = 1, \dots, d$ are not asymptotically uncorrelated.
\end{remark}

\section*{Complete Measure of Mutual Dependence Using Weight Function $w_2$}

Except that $\mathcal{U}_n(\mathbf{X})$ requires the additional $d$-th moment condition $\textrm{E}|X_1  \dots  X_d| < \infty$ to be simplified, $\mathcal{U}(X)$ is in an extremely complicated form. Even when $d = 3$, $\mathcal{U}(X)$ already has 12 different terms as follows.
\begin{eqnarray*}
  \mathcal{U}(X) &=& \|\phi_{X}(t) - \phi_{\widetilde{X}}(t)\|^2_{w_2} \\
  &=& - \textrm{E}|X_1 - X_1'||X_2 - X_2'||X_3 - X_3'| \\
  && + \, 2\textrm{E}|X_1 - X_1'||X_2 - X_2''||X_3 - X_3'''| \\
  && - \, \textrm{E}|X_1 - X_1'|\textrm{E}|X_2 - X_2'|\textrm{E}|X_3 - X_3'| \\
  && + \, \textrm{E}|X_1 - X_1'||X_2 - X_2'| + \, \textrm{E}|X_1 - X_1'||X_3 - X_3'| + \, \textrm{E}|X_2 - X_2'||X_3 - X_3'| \\
  && - \, 2\textrm{E}|X_1 - X_1'||X_2 - X_2''| - \, 2\textrm{E}|X_1 - X_1'||X_3 - X_3'''| - \, 2\textrm{E}|X_2 - X_2''||X_3 - X_3'''| \\
  && + \, \textrm{E}|X_1 - X_1'|\textrm{E}|X_2 - X_2'| + \, \textrm{E}|X_1 - X_1'|\textrm{E}|X_3 - X_3'| + \, \textrm{E}|X_2 - X_2'|\textrm{E}|X_3 - X_3'| \\
  &=& - \textrm{E}|X_1 - X_1'||X_2 - X_2'||X_3 - X_3'| \\
  && + \, 2\textrm{E}|X_1 - X_1'||X_2 - X_2''||X_3 - X_3'''| \\
  && - \, \textrm{E}|X_1 - X_1'|\textrm{E}|X_2 - X_2'|\textrm{E}|X_3 - X_3'| \\
  && + \, \sum_{1 \leq i < j \leq 3} \textrm{E}|X_i - X_i'||X_j - X_j'| \\
  && - \, 2\sum_{1 \leq i < j \leq 3} \textrm{E}|X_i - X_i'||X_j - X_j''| \\
  && + \, \sum_{1 \leq i < j \leq 3} \textrm{E}|X_i - X_i'|\textrm{E}|X_j - X_j'|.
\end{eqnarray*}

In general, the number of different terms in $\mathcal{U}(X)$ grows exponentially as $d$ increases. Basically, we will see all combinations of all components in all moments as expectations.

\end{document}